\documentclass[12pt,leqno]{article}
\usepackage{amssymb}
\usepackage[mathscr]{eucal}
\usepackage{amsmath,amssymb,latexsym,theorem,bbm}
\usepackage{color}
\usepackage{appendix}

\newcommand{\SC}{\scriptstyle}

\newcommand{\CC}{\mathbb{C}}

\newcommand{\NN}{\mathbb{N}}
\newcommand{\RR}{\mathbb{R}}

\newcommand{\ZZ}{\mathbb{Z}}

\newcommand{\tm}{\widetilde{m}}

\newcommand{\hvartheta}{{\widehat{\vartheta}}}

\newcommand{\cB}{{\mathcal B}}

\newcommand{\cD}{{\mathcal D}}
\newcommand{\cE}{{\mathcal E}}

\newcommand{\cY}{{\mathcal Y}}
\newcommand{\cZ}{{\mathcal Z}}
\newcommand{\cW}{{\mathcal W}}

\newcommand{\dd}{\mathrm{d}}
\newcommand{\ee}{\mathrm{e}}

\newcommand{\ii}{\mathrm{i}}

\newcommand{\EE}{\operatorname{\mathbb{E}}}
\newcommand{\PP}{{\operatorname{\mathbb{P}}}}

\newcommand{\oo}{\operatorname{o}}
\newcommand{\OO}{\operatorname{O}}

\renewcommand{\Re}{\operatorname{Re}}
\renewcommand{\Im}{\operatorname{Im}}
\newcommand{\Res}{\operatorname{Res}}

\newcommand{\tY}{\widetilde{Y}}

\renewcommand{\leq}{\leqslant}
\renewcommand{\geq}{\geqslant}

\newcommand{\stoch}{\stackrel{\PP}{\longrightarrow}}
\newcommand{\distr}{\stackrel{\cD}{\longrightarrow}}

\newcommand{\as}{\stackrel{{\mathrm{a.s.}}}{\longrightarrow}}

\newcommand{\bbone}{\mathbbm{1}}

\newcommand{\proofend}{\hfill\mbox{$\Box$}}

\numberwithin{equation}{section}

\theoremstyle{change} \theorembodyfont{\em}
\newtheorem{Lem}{Lemma.}[section]
\newtheorem{Thm}[Lem]{Theorem.}

\newtheorem{Cor}[Lem]{Corollary.}

\theorembodyfont{\rm}
\newtheorem{Rem}[Lem]{Remark.}

\begin{document}

\begin{center}
 {\bfseries\Large One-parameter statistical model for linear stochastic}\\[2mm]
 {\bfseries\Large differential equation with time delay} \\[5mm]
 {\sc\large J\'anos Marcell $\text{Benke}^*$ \ and \ Gyula Pap}
\end{center}

\vskip0.2cm

\noindent
 Bolyai Institute, University of Szeged,
 Aradi v\'ertan\'uk tere 1, H--6720 Szeged, Hungary.

\noindent e--mails: jbenke@math.u-szeged.hu (J. M. Benke),
                    papgy@math.u-szeged.hu (G. Pap).

\noindent * Corresponding author.

\renewcommand{\thefootnote}{}
\footnote{\textit{2010 Mathematics Subject Classifications\/}:
          62B15, 62F12.}
\footnote{\textit{Key words and phrases\/}:
 likelihood function; local asymptotic normality; local asymptotic mixed
 normality; periodic local asymptotic mixed normality; local asymptotic
 quadraticity; maximum likelihood estimator; stochastic differential equations;
 time delay.}

\vspace*{-7mm}


\begin{abstract}
Assume that we observe a stochastic process \ $(X(t))_{t\in[-r,T]}$, \ which
 satisfies the linear stochastic delay differential equation
 \[
   \dd X(t)
   = \vartheta \int_{[-r,0]} X(t + u) \, a(\dd u) \, \dd t + \dd W(t) ,
   \qquad t \geq 0 ,
 \]
 where \ $a$ \ is a finite signed measure on \ $[-r, 0]$.
\ The local asymptotic properties of the likelihood function are studied.
Local asymptotic normality is proved in case of \ $v_\vartheta^* < 0$, \ local
 asymptotic quadraticity is shown if \ $v_\vartheta^* = 0$, \ and, under some additional
 conditions, local asymptotic mixed normality or periodic local asymptotic mixed
 normality is valid if \ $v_\vartheta^* > 0$, \ where \ $v_\vartheta^*$ \ is an
 appropriately defined quantity.
As an application, the asymptotic behaviour of the maximum likelihood estimator
 \ $\hvartheta_T$ \ of \ $\vartheta$ \ based on \ $(X(t))_{t\in[-r,T]}$ \ can be derived
 \ as \ $T \to \infty$.
\end{abstract}

\section{Introduction}

Consider the linear stochastic delay differential equation (SDDE)
 \begin{equation}\label{SDDE}
  \begin{cases}
   \dd X(t)
   = \vartheta \int_{[-r,0]} X(t + u) \, a(\dd u) \, \dd t
     + \dd W(t) ,
    & t \in \RR_+ , \\
   X(t) = X_0(t) , & t \in [-r, 0] , 
  \end{cases}
 \end{equation}
 where \ $r \in (0, \infty)$, \ $(W(t))_{t\in\RR_+}$ \ is a standard Wiener
 process, \ $\vartheta \in \RR$, \ and \ $a$ \ is a finite signed measure on
 \ $[-r, 0]$ \ with \ $a \ne 0$, \ and \ $(X_0(t))_{t\in[-r,0]}$ \ is a
 continuous process independent of \ $(W(t))_{t\in\RR_+}$. 
\ The SDDE \eqref{SDDE} can also be written in the integral form
 \begin{equation}\label{iSDDE}
  \begin{cases}
   X(t) = X_0(0)
          + \vartheta \int_0^t \int_{[-r,0]} X(s + u) \, a(\dd u) \, \dd s
          + W(t) ,
    & t \in \RR_+ , \\
   X(t) = X_0(t) , & t \in [-r, 0] . 
  \end{cases}
 \end{equation}
Equation \eqref{SDDE} is a special case of the affine stochastic delay differential
 equation
 \begin{equation}\label{ASDDE}
  \begin{cases}
   \dd X(t) = \int_{-r}^0 X(t + u) \, a_\vartheta(\dd u) \, \dd t + \dd W(t) ,
    & t \in \RR_+ , \\
   X(t) = X_0(t) , & t \in [-r, 0] , 
  \end{cases}
 \end{equation}
 where \ $r > 0$, \ and \ for each \ $\vartheta \in \Theta$, \ $a_\vartheta$ \ is a
 finite signed measure on \ $[-r, 0]$, \ see Gushchin and K\"uchler \cite{GusKuc2003}.
In that paper local asymptotic normality (LAN) has been proved for stationary solutions.
In Gushchin and K\"uchler \cite{GusKuc1999}, the special case of \eqref{ASDDE} has
 been studied with \ $r = 1$, \ $\Theta = \RR^2$, \ and
 \ $a_\vartheta = \vartheta_1 \delta_0 + \vartheta_2 \delta_{-1}$ \ for
 \ $\vartheta = (\vartheta_1, \vartheta_2)$, \ where \ $\delta_x$ \ denotes the Dirac
 measure concentrated at \ $x \in \RR$, \ and they described the local properties of the
 likelihood function for the whole parameter space \ $\RR^2$. 
\ In Benke and Pap \cite{BenPap2015}, a special case has been studied, where
 \ $r = 1$ \ and \ $a_\vartheta$ \ is the Lebesgue measure multiplied by
 \ $\vartheta \in \RR$.
 
In each of the above papers, LAN has been proved in case of \ $v_0(\vartheta) < 0$,
 \ where \ $v_0(\vartheta)$ \ is the real part of the right most characteristic roots of
 the corresponding deterministic homogeneous delay differential equation, see \eqref{v0}. 
It turns out that in case of equation \eqref{SDDE}, LAN holds whenever
 \ $v_\vartheta^* < 0$, \ where \ $v_\vartheta^*$ \ is defined in \eqref{v*m*},
 see Theorem \eqref{LAN}, but it can happen that \ $v_0(\vartheta) = 0$, \ see the
 example in Remark \ref{v*}.
Moreover, local asymptotic quadraticity (LAQ) is shown if \ $v_\vartheta^* = 0$, \ and,
 under some additional conditions, local asymptotic mixed normality (LAMN) or periodic
 local asymptotic mixed normality (PLAMN) is valid if \ $v_\vartheta^* > 0$, \ see
 Theorems \ref{LAQ} and \ref{PLAMN}.
Note that in Theorems \ref{LAQ} and \ref{PLAMN} we have
 \ $v_\vartheta^* = v_0(\vartheta)$, see Remark \ref{v*}.
The definition of LAN, LAQ, LAMN and PLAMN can be found in Le Cam and Yang
 \cite{LeCamYang} and Gushchin and K\"uchler \cite{GusKuc1999}.
 
The solution \ $(X^{(\vartheta)}(t))_{t\in\RR_+}$ \ of \eqref{SDDE} exists, is
 pathwise uniquely determined and can be represented as
 \begin{equation}\label{solution}
  \begin{aligned}
   X^{(\vartheta)}(t)
   &= x_{0,\vartheta}(t) X_0(0)
      + \vartheta 
        \int_{[-r,0]}
         \int_u^0 x_{0,\vartheta}(t + u - s) X_0(s) \, \dd s \, a(\dd u) \\
   &\quad
      + \int_{[0,t]} W(t - s) \, \dd x_{0,\vartheta}(s) ,
    \qquad t \in \RR_+ ,
  \end{aligned}
 \end{equation}
 where \ $(x_{0,\vartheta}(t))_{t\in[-r,\infty)}$ \ denotes the so-called
 fundamental solution of the deterministic homogeneous delay differential
 equation
 \begin{equation}\label{fund}
  \begin{cases}
   x(t) = x_0(0)
          + \vartheta \int_0^t \int_{[-r,0]} x(s + u) \, a(\dd u) \, \dd s ,
    & t \in \RR_+ , \\
   x(t) = x_0(t) , & t \in [-r, 0] . 
  \end{cases}
 \end{equation}
 with initial function
 \[
   x_0(t) := \begin{cases}
              0 , & t \in [-r, 0) , \\
              1 , & t = 0 ,
             \end{cases}
 \]
 which means that \ $x_{0,\vartheta}$ \ is absolutely continuous on \ $\RR_+$,
 \ $x_{0,\vartheta}(t) = 0$ \ for \ $t \in [-r, 0)$,
 \ $x_{0,\vartheta}(0) = 1$, \ and
 \ $\dot{x}_{0,\vartheta}(t)
    = \vartheta \int_{[-r,0]} x_{0,\vartheta}(t + u) \, a(\dd u)$
 \ for Lebesgue-almost all \ $t \in \RR_+$.
\ The domain of integration in the last integral in \eqref{solution} includes
 zero, i.e.,
 \[
   \int_{[0,t]} W(t - s) \, \dd x_{0,\vartheta}(s)
   = W(t) + \int_{(0,t]} W(t - s) \, \dd x_{0,\vartheta}(s)
   = \int_0^t x_{0,\vartheta}(t - s) \, \dd W(s) , \qquad t \in \RR_+ .
 \]
In the trivial case of \ $\vartheta = 0$, \ we have \ $x_{0,0}(t) = 1$ \ for
 all \ $t \in \RR_+$, \ and \ $X^{(0)}(t) = X_0(0) + W(t)$ \ for all
 \ $t \in \RR_+$.
\ The asymptotic behaviour of \ $x_{0,\vartheta}(t)$ \ as \ $t \to \infty$ \ is
 connected with the so-called characteristic function
 \ $h_\vartheta : \CC \to \CC$, \ given by
 \begin{equation}\label{char1}
  h_\vartheta(\lambda)
  := \lambda - \vartheta \int_{[-r,0]} \ee^{\lambda u} \, a(\dd u) ,
  \qquad \lambda \in \CC ,
 \end{equation}
 and the set \ $\Lambda_\vartheta$ \ of the (complex) solutions of the 
 so-called characteristic equation for \eqref{fund},
 \begin{equation}\label{char2}
  \lambda - \vartheta \int_{[-r,0]} \ee^{\lambda u} \, a(\dd u) = 0 .
 \end{equation} 
Note that a complex number \ $\lambda$ \ solves \eqref{char2} if and
 only if \ $(\ee^{\lambda t})_{t\in[-r,\infty)}$ \ solves \eqref{fund}
 with initial function \ $x_0(t) = \ee^{\lambda t}$, \ $t \in [-r, 0]$.
\ We have \ $\Lambda_\vartheta \ne \emptyset$,
 \ $\overline{\Lambda_\vartheta} = \Lambda_\vartheta$, \ and 
 \ $\Lambda_\vartheta$ \ consists of isolated points.
Moreover, \ $\Lambda_\vartheta$ \ is countably infinite except the case where
 \ $a$ \ is concentrated at 0, \ or \ $\vartheta = 0$.
\ Further, for each \ $c \in \RR$, \ the set
 \ $\{\lambda \in \Lambda_\vartheta : \Re(\lambda) \geq c\}$ \ is finite.
In particular,
 \begin{equation}\label{v0}
   v_0(\vartheta)
   := \sup\{\Re(\lambda) : \lambda \in \Lambda_\vartheta\} < \infty .
 \end{equation}
For \ $\lambda \in \Lambda_\vartheta$, \ denote by \ $m_\vartheta(\lambda)$
 \ the multiplicity of \ $\lambda$ \ as a solution of \eqref{char2}.
 
The Laplace transform of \ $(x_{0,\vartheta}(t))_{t\in\RR_+}$ \ is given by
 \[
   \int_0^\infty \ee^{-\lambda t} x_{0,\vartheta}(t) \, \dd t
   = \frac{1}{h_\vartheta(\lambda)} , \qquad \lambda \in \CC , \qquad
   \Re(\lambda) > v_0(\vartheta) .   
 \]
Based on the inverse Laplace transform and Cauchy's residue theorem,
 the following crucial lemma can be shown (see, e.g., Diekmann et al.\
 \cite[Lemma 5.1 and Theorem 5.4]{DvGVLW} or Gushchin and K\"uchler
 \cite[Lemma 2.1]{GusKuc2000}).  
 
\begin{Lem}\label{res}
For each \ $\vartheta \in \RR$ \ and each \ $c \in \RR$, \ the fundamental
 solution \ $(x_{0,\vartheta}(t))_{t\in[-r,\infty)}$ \ of \eqref{fund} can be
 represented in the form
 \[
   x_{0,\vartheta}(t)
   =\sum_{\underset{\SC\Re(\lambda)\geq c}
                   {\lambda \in \Lambda_\vartheta}}
     \underset{\SC z=\lambda}{\Res}\left(\frac{\ee^{zt}}{h_\vartheta(z)}\right)
    + \psi_{\vartheta,c}(t)
   =\sum_{\underset{\SC\Re(\lambda)\geq c}
                   {\lambda \in \Lambda_\vartheta}}
     p_{\vartheta,\lambda}(t) \, \ee^{\lambda t}
    + \psi_{\vartheta,c}(t) , \qquad \text{as \ $t \to \infty$,}
 \]
 where \ $\psi_{\vartheta,c} : \RR_+ \to \RR$ \ is a continuous function with
 \ $\psi_{\vartheta,c}(t) = \oo(\ee^{c t})$ \ as \ $t \to \infty$, \ and for
 each \ $\vartheta \in \RR$ \ and each \ $\lambda \in \Lambda_\vartheta$,
 \ $p_{\vartheta,\lambda}$ \ is a complex-valued polynomial of degree
 \ $m_\vartheta(\lambda) - 1$ \ with
 \ $p_{\vartheta,\overline{\lambda}} = \overline{p_{\vartheta,\lambda}}$.
\ More exactly,
 \[
   p_{\vartheta,\lambda}(t) 
   = \sum_{\ell=0}^{m_\vartheta(\lambda)-1}
      \frac{A_{\vartheta,-1-\ell}(\lambda)}{\ell!} \, t^\ell ,
 \]
 where \ $A_{\vartheta,k}(\lambda)$,
 \ $k \in \{-m_\vartheta(\lambda), -m_\vartheta(\lambda)+1, \ldots\}$ \ denotes
 the coefficients of the Laurent's series of \ $1/h_\vartheta(z)$ \ at
 \ $z = \lambda$, \ i.e.,
 \[
   \frac{1}{h_\vartheta(z)}
   = \sum_{k=-m_\vartheta(\lambda)}^\infty
      A_{\vartheta,k}(\lambda) (z - \lambda)^k 
 \]
 in a neighborhood of \ $\lambda$.
\end{Lem}

As a consequence, for any \ $c > v_0(\vartheta)$, \ we have
 \ $x_{0,\vartheta}(t) = \OO(\ee^{c t})$, \ as \ $t \to \infty$.
\ In particular, \ $(x_{0,\vartheta}(t))_{t\in\RR_+}$ \ is square integrable
 if (and only if, see Gushchin and K\"uchler \cite{GusKuc2000})
 \ $v_0(\vartheta) < 0$.

\section{Radon--Nikodym derivatives}
\label{section_RN}

From this section, we will consider the SDDE \eqref{SDDE} with fixed
 continuous initial process \ $(X_0(t))_{t\in[-r,0]}$.
\ Further, for all \ $T \in \RR_{++}$, \ let \ $\PP_{\vartheta,T}$ \ be the
 probability measure induced by \ $(X^{(\vartheta)}(t))_{t\in[-r,T]}$ \ on
 \ $(C([-r, T]), \cB(C([-r, T])))$. 
\ In order to calculate Radon--Nikodym derivatives
 \ $\frac{\dd \PP_{\theta,T}}{\dd \PP_{\vartheta,T}}$ \ for certain
 \ $\theta, \vartheta \in \RR$, \ we need the following statement, which
 can be derived from formula (7.139) in Section 7.6.4 of Liptser and Shiryaev
 \cite{LipShiI}.

\begin{Lem}\label{RN}
Let \ $\theta, \vartheta \in \RR$.
\ Then for all \ $T \in \RR_{++}$, \ the measures \ $\PP_{\theta,T}$ \ and
 \ $\PP_{\vartheta,T}$ \ are absolutely continuous with respect to each other,
 and
 \begin{align*}
  \log \frac{\dd \PP_{\theta,T}}{\dd \PP_{\vartheta,T}}
   (X^{(\vartheta)}|_{[-r,T]})
  &= (\theta - \vartheta) \int_0^T Y^{(\vartheta)}(t) \, \dd X^{(\vartheta)}(t)
     - \frac{1}{2} (\theta^2 - \vartheta^2)
       \int_0^T Y^{(\vartheta)}(t)^2 \, \dd t \\
  &= (\theta - \vartheta) \int_0^T Y^{(\vartheta)}(t) \, \dd W(t)
     - \frac{1}{2}
       (\theta - \vartheta)^2
       \int_0^T Y^{(\vartheta)}(t)^2 \, \dd t
 \end{align*}
 with
 \[
   Y^{(\vartheta)}(t) := \int_{[-r,0]} X^{(\vartheta)}(t+u) \, a(\dd u) ,
   \qquad t \in \RR_+ .
 \]
\end{Lem}

In order to investigate local asymptotic properties of the family
 \begin{equation}\label{cET}
  (\cE_T)_{T\in\RR_{++}}
  := \big( C(\RR_+), \cB(C(\RR_+)) ,
           \{\PP_{\vartheta,T} : \vartheta \in \RR\} \big)_{T\in\RR_{++}}
 \end{equation}
 of statistical experiments, we derive the following corollary.

\begin{Cor}\label{RN_Cor}
For each \ $\vartheta \in \RR$, \ $T \in \RR_{++}$,
 \ $r_{\vartheta,T} \in \RR$ \ and \ $h_T \in \RR$, \ we have
 \[
   \log \frac{\dd \PP_{\vartheta+r^{(\vartheta,T)}h_T,T}}
             {\dd \PP_{\vartheta,T}}(X^{(\vartheta)}|_{[-r,T]})
   = h_T \Delta_{\vartheta,T} - \frac{1}{2} h_T^2 J_{\vartheta,T} ,
 \]
 with
 \[
   \Delta_{\vartheta,T}
   := r_{\vartheta,T} \int_0^T Y^{(\vartheta)}(t) \, \dd W(t) , \qquad
   J_{\vartheta,T}
   := r_{\vartheta,T}^2 \int_0^T Y^{(\vartheta)}(t)^2 \, \dd t .
 \]
\end{Cor}

\section{Local asymptotics of likelihood ratios}
\label{section_LALR}

For each \ $\lambda \in \Lambda_\vartheta$, \ denote by
 \ $\tm_\vartheta(\lambda)$ \ the degree of the complex-valued polynomial
 \[
   P_{\vartheta,\lambda}(t)
   := \sum_{\ell=0}^{m_\vartheta(\lambda)-1}
       c_{\vartheta,\lambda,\ell} \, t^\ell
 \]
 with
 \[
   c_{\vartheta,\lambda,\ell}
   := \frac{1}{\ell!}
      \int_{[-r, 0]}
       \underset{\SC z=\lambda}{\Res}
        \biggl(\frac{(z-\lambda)^\ell \ee^{zu}}{h_\vartheta(z)}\biggr)
       a(\dd u)
   = \frac{1}{\ell!}
     \sum_{j=0}^{m_\vartheta(\lambda)-1-\ell}
      \frac{A_{\vartheta,-j-1-\ell}(\lambda)}{j!}
      \int_{[-r, 0]} u^j \ee^{\lambda u} \, a(\dd u) ,
 \]
 where the degree of the zero polynomial is defined to be \ $-\infty$.
\ Put
 \begin{equation}\label{v*m*}
   v_\vartheta^*
   := \sup\{\Re(\lambda) : \lambda \in \Lambda_\vartheta , \;
                           \tm_\vartheta(\lambda) \geq 0\} , \qquad
   m_\vartheta^*
   := \max\{\tm_\vartheta(\lambda) : \lambda \in \Lambda_\vartheta, \;
                                     \Re(\lambda) = v_\vartheta^*\} ,
 \end{equation}
 where \ $\sup \emptyset := -\infty$ \ and \ $\max \emptyset := -\infty$.
 
\begin{Thm}\label{LAN}
If \ $\vartheta \in \RR$ \ with \ $v_\vartheta^* < 0$, \ then the family
 \ $(\cE_T)_{T\in\RR_{++}}$ \ of statistical experiments given in \eqref{cET}
 is LAN at \ $\vartheta$ \ with scaling \ $r_{\vartheta,T} = T^{-1/2}$,
 \ $T \in \RR_{++}$, \ and with
 \[
   J_\vartheta
   = \int_0^\infty
      \left(\int_{[-r,0]} x_{0,\vartheta}(t + u) \, a(\dd u)\right)^2 \dd t . 
 \]
Particularly, if \ $a([-r, 0]) = 0$, \ then \ $v_0^* = -\infty$,
 \ $m_0^* = -\infty$, \ and the family \ $(\cE_T)_{T\in\RR_{++}}$ \ of
 statistical experiments given in \eqref{cET} is LAN at \ $0$ \ with scaling
 \ $r_{0,T} = T^{-1/2}$, \ $T \in \RR_{++}$, \ and with
 \[
   J_0 = \int_0^r a([-t, 0])^2 \, \dd t . 
 \]
\end{Thm}

\begin{Thm}\label{LAQ}
If \ $\vartheta \in \RR$ \ with \ $v_\vartheta^* = 0$, \ then the family
 \ $(\cE_T)_{T\in\RR_{++}}$ \ of statistical experiments given in \eqref{cET}
 is LAQ at \ $\vartheta$ \ with scaling
 \ $r_{\vartheta,T} = T^{-m_\vartheta^*-1}$ \ and with
 \begin{gather*}
  \Delta_\vartheta
  = \sum_{\underset{\SC\tm_\vartheta(\lambda)=m_\vartheta^*}
                   {\lambda \in \Lambda_\vartheta\cap(\ii\RR)}}
     c_{\vartheta,\lambda,m_\vartheta^*}
     \int_0^1
      \cZ_{\Im(\lambda),m_\vartheta^*}(s)
      \, \dd \overline{\cZ_{\Im(\lambda),0}(s)} , \\
  J_\vartheta
  = \sum_{\underset{\SC\tm_\vartheta(\lambda)=m_\vartheta^*}
                   {\lambda \in \Lambda_\vartheta\cap(\ii\RR)}}
     |c_{\vartheta,\lambda,m_\vartheta^*}|^2
     \int_0^1 |\cZ_{\Im(\lambda),m_\vartheta^*}(s)|^2 \, \dd s ,
 \end{gather*}
 with
 \[
   \cZ_{\varphi,0}
   := \begin{cases}
       \cW, & \text{if \ $\varphi = 0$,} \\
       \frac{1}{\sqrt{2}}
       \bigl(\cW_{\varphi,\Re} + \ii \cW_{\varphi,\Im}\bigr),
        & \text{if \ $\varphi \in \RR_{++}$,} \\[1mm]
       \overline{\cZ_{-\varphi,0}},
        & \text{if \ $\varphi \in \RR_{--}$,}
      \end{cases}
 \]
 where \ $(\cW(s))_{s\in[0,1]}$, \ $(\cW_{\varphi,\Re}(s))_{s\in[0,1]}$
 \ and \ $(\cW_{\varphi,\Im}(s))_{s\in[0,1]}$, \ $\varphi \in \RR_{++}$,
 \ are independent standard Wiener processes, and
 \[
   \cZ_{\varphi,\ell}(s)
   := \int_0^s (s - u)^\ell \, \dd \cZ_{\varphi,0}(u) , \qquad
   s \in [0, 1] , \quad \varphi \in \RR , \quad \ell \in \NN .
 \]
Particularly, if \ $a([-r, 0]) \ne 0$, \ then \ $v_0^* = 0$, \ $m_0^* = 0$,
 \ and the family \ $(\cE_T)_{T\in\RR_{++}}$ \ of statistical experiments given
 in \eqref{cET} is LAQ at \ $0$ \ with scaling \ $r_{0,T} = T^{-1}$,
 \ $T \in \RR_{++}$, \ and with 
 \[
   \Delta_0 = a([-r, 0]) \int_0^1 \cW(s) \, \dd \cW(s) , \qquad
   J_0 = a([-r, 0])^2 \int_0^1 \cW(s)^2 \, \dd s . 
 \]
\end{Thm}

Note that \ $\Delta_\vartheta$ \ is almost surely real-valued, since
 \[
   c_{\vartheta,\overline{\lambda},m_\vartheta^*}
     \int_0^1
      \cZ_{\Im(\overline{\lambda}),m_\vartheta^*}(s)
      \, \dd \overline{\cZ_{\Im(\overline{\lambda}),0}(s)}
   = \overline{c_{\vartheta,\lambda,m_\vartheta^*}}
               \int_0^1
                \overline{\cZ_{\Im(\lambda),m_\vartheta^*}(s)}
                \, \dd \cZ_{\Im(\lambda),0}(s) , \qquad
   \lambda \in \Lambda_\vartheta .
 \]

\begin{Thm}\label{PLAMN}
Let \ $\vartheta \in \RR$ \ with \ $v_\vartheta^* > 0$.
\ If
 \[
   H_\vartheta
   := \{\Im(\lambda)
        : \lambda \in \Lambda_\vartheta \cap (v_\vartheta^* + \ii \RR_{++}) , \;
          \tm_\vartheta(\lambda) = m_\vartheta^*\}
   \ne \emptyset ,
 \]
 and the numbers in \ $H_\vartheta$ \ have a common divisor \ $D_\vartheta$,
 \ then the family \ $(\cE_T)_{T\in\RR_{++}}$
 \ of statistical experiments given in \eqref{cET} is PLAMN at \ $\vartheta$ \ with
 period $\frac{2 \pi}{D_\vartheta}$,  
 \ with scaling \ $r_{\vartheta,T} = T^{-m_\vartheta^*} \ee^{-v_\vartheta^* T}$,
 \ $T \in \RR_{++}$, \ and with
 \[
   \Delta_\vartheta(d) = Z \sqrt{J_\vartheta(d)} , \qquad
   J_\vartheta(d)
   = \int_0^\infty
      \ee^{-2v_\vartheta^* t}
      \Re\Biggl(\sum_{\underset{\SC\tm_\vartheta(\lambda)=m_\vartheta^*}
                    {\lambda\in\Lambda_\vartheta\cap(v_\vartheta^*+\ii\RR)}}
                c_{\vartheta,\lambda,m_\vartheta^*}
                U^{(\vartheta)}_\lambda
                \ee^{\ii(d-t)\Im(\lambda)}\Biggr)^2
      \dd t ,
 \]
 for \ $d \in [0, \frac{2 \pi}{D_\vartheta})$, \ where
 \[
   U^{(\vartheta)}_\lambda
   := X_0(0)  
      + v_\vartheta^*
        \int_{[-r,0]} \int_u^0 \ee^{-\lambda(s-u)} X_0(s) \, \dd s \, a(\dd u)
      + \int_0^\infty \ee^{-\lambda s} \, \dd W(s) , \qquad \lambda \in \CC ,
 \]
 and \ $Z$ \ is a standard normally distributed random variable independent
 of \ $(X_0(t))_{t\in[-r,0]}$ \ and \ $(W(t))_{t\in\RR_+}$.

If \ $H_\vartheta = \emptyset$, \ then the family \ $(\cE_T)_{T\in\RR_{++}}$
 \ of statistical experiments given in \eqref{cET} is LAMN at \ $\vartheta$
 \ with scaling \ $r_{\vartheta,T} = T^{-m_\vartheta^*} \ee^{-v_\vartheta^* T}$,
 \ $T \in \RR_{++}$, \ and with
 \[
   \Delta_\vartheta = Z \sqrt{J_\vartheta} , \qquad
   J_\vartheta
   = \frac{c_{\vartheta,v_\vartheta^*,m_\vartheta^*}^2}{2v_\vartheta^*}
     (U^{(\vartheta)}_{v_\vartheta^*})^2 .
 \]
\end{Thm}

\begin{Rem}\label{v*}
According to the definition of \ $v_0(\vartheta)$ \ and \ $v_\vartheta^*$, \ we
 obtain \ $v_0(\vartheta) \ge v_\vartheta^*$. 
\ The aim of the following discussion is to show that
 \ $v_0(\vartheta) > v_\vartheta^*$ \ if and only if
 \ $\{\lambda \in \Lambda_\vartheta : \Re(\lambda) \geq 0\} = \{0\}$ \ and
 \ $P_{\vartheta,0} = 0$ \ (and hence \ $v_0(\vartheta) = 0$ \ and
 \ $v_\vartheta^* < 0$).
\ Indeed, if \ $v_0(\vartheta) > v_\vartheta^*$ \ then for all
 \ $\lambda_0 \in \Lambda_\vartheta$ \ with \ $\Re(\lambda_0) > v_\vartheta^*$
 \ we have \ $P_{\vartheta,\lambda_0} = 0$, \ implying
 \[
   c_{\vartheta ,\lambda_0 , m_\vartheta(\lambda_0) -1}
   = \frac{1}{(m_\vartheta(\lambda_0) -1)!}
     A_{\vartheta,-m_\vartheta(\lambda_0)}(\lambda_0)
     \int_{[-r, 0]} \ee^{\lambda_0 u} \, a(\dd u)
   = 0 .
 \]
Here \ $A_{\vartheta,-m_\vartheta(\lambda_0)}(\lambda) \ne 0$, \ since
 \ it is the leading coefficient of the polynomial \ $p_{\vartheta,\lambda_0}$ \ of
 degree \ $m_\vartheta(\lambda_0) -1$, \ hence
 \ $c_{\vartheta ,\lambda_0 , m_\vartheta(\lambda_0) -1} = 0$ \ yields
 \ $\int_{[-r,0]} e^{\lambda_0 u} a( \dd u) = 0$.
\  Taking into account of the characteristic equation, we get \ $\lambda_0 = 0$,
 \ hence
 \ $\{\lambda \in \Lambda_\vartheta : \Re(\lambda) \geq v_\vartheta^*\} = \{0\}$. 
\ Clearly, this yields also \ $v_0(\vartheta) = 0$ \ and
 \ $v_\vartheta^* < 0$, \ and hence
 \ $\{\lambda \in \Lambda_\vartheta : \Re(\lambda) \geq 0\} = \{0\}$.
\ Conversely, if \ $\{\lambda \in \Lambda_\vartheta : \Re(\lambda) \geq 0\} = \{0\}$
 \ and \ $P_{\vartheta,0} = 0$, \ then, by definition, \ $v_0(\vartheta) = 0$ \ and
 \ $v_\vartheta^* < 0$.

In particular, if \ $\{\lambda \in \Lambda_\vartheta : \Re(\lambda) \geq 0\} = \{0\}$
 \ and \ $m_\vartheta(0) = 1$, \ then \ $v_0(\vartheta) > v_\vartheta^*$ \ is equivalent
 to \ $a([-r, 0]) = 0$, \ since 
 \ $P_{\vartheta,0} = c_{\vartheta ,0 , 0}
    = \int_{[-r, 0]} \ee^{\lambda_0 u} \, a(\dd u) = a([-r, 0])$.
 
An example for this situation is, when \ $r = 2 \pi$ , \ $a(\dd u) = \sin(u) \dd u$
 \ and \ $\vartheta \in \big( 0,\frac{1}{\pi}\big)$.
\ This can be derived, applying usual methods (e.g., argument principle in complex
 analysis and the existence of local inverses of holomorphic functions), see,
 e.g., Rei{\ss} \cite{Rei}.
\end{Rem}

\begin{Rem}
Using these results, we can give the asymptotic behaviour of the maximum likelihood
 estimator of \ $\vartheta$ \ based on the observation \ $\big(X(t)\big)_{t \in [-1,T]}$
 \ for some fixed \ $ T > 0$, \ see, e.g., Benke and Pap \cite[Section 5]{BenPap2015}.
\end{Rem}

\section{Proofs}
\label{section_proofs}

For each \ $\vartheta \in \RR$ \ and each \ $t \in [r, \infty)$, \ by
 \eqref{solution}, we have
 \begin{align*}
  Y^{(\vartheta)}(t)
  &= X_0(0) \int_{[-r,0]} x_{0,\vartheta}(t + u) \, \dd u
     + \int_{[-r,0]} \int_{[0,t+u]}
        W(t + u - s) \, \dd x_{0,\vartheta}(s) \, a(\dd u) \\
  &\quad
     + \vartheta \int_{[-r,0]} \int_{[-r,0]} \int_v^0
        x_{0,\vartheta}(t + u + v - s) X_0(s)
        \, \dd s \, a(\dd v) \, a(\dd u) .
 \end{align*}
Here we have
 \begin{align*}
  &\int_{[-r,0]} \int_{[-r,0]} \int_v^0
    x_{0,\vartheta}(t + u + v - s) X_0(s)
    \, \dd s \, a(\dd v) \, a(\dd u) \\
  &\qquad\qquad
   = \int_{[-r,0]} \int_{[-r,0]} \int_v^0
      x_{0,\vartheta}(t + u + v - s) X_0(s)
      \, \dd s \, a(\dd u) \, a(\dd v) \\
  &\qquad\qquad
   = \int_{[-r,0]} \int_v^0 X_0(s) \int_{[-r,0]}
      x_{0,\vartheta}(t + u + v - s) \, a(\dd u) \, \dd s \, a(\dd v) ,
 \end{align*}
 and
 \begin{align*}
  &\int_{[-r,0]} \int_0^{t+u}
    x_{0,\vartheta}(t + u - s) \, \dd W(s) \, a(\dd u) \\
  &= \int_0^{t-r} \int_{[-r,0]}
      x_{0,\vartheta}(t + u - s) \, a(\dd u) \, \dd W(s)
     + \int_{t-r}^t \int_{[s-t,0]}
        x_{0,\vartheta}(t + u - s) \, a(\dd u) \, \dd W(s) \\
  &= \int_0^t \int_{[-r,0]}
      x_{0,\vartheta}(t + u - s) \, a(\dd u) \, \dd W(s) ,
 \end{align*}
 since \ $t \in [r, \infty)$, \ $s \in [t-r, t]$ \ and \ $u \in [-r, s-t)$
 \ imply \ $t + u - s \in [-r, 0)$, \ and hence
 \ $x_{0,\vartheta}(t + u - s) = 0$.
\ Consequently, the process \ $\bigl(Y^{(\vartheta)}(t)\bigr)_{t\in[r,\infty)}$
 \ has a representation
 \begin{equation}\label{Yvartheta}
  \begin{aligned}
   Y^{(\vartheta)}(t)
   &= y_\vartheta(t) X_0(0)
      + \vartheta
        \int_{[-r,0]} \int_v^0
         y_\vartheta(t + v - s) X_0(s) \, \dd s \, a(\dd v)
      + \int_0^t y_\vartheta(t - s) \, \dd W(s)
  \end{aligned}
 \end{equation}
 for \ $t \in [r,\infty)$ \ with
 \[
   y_\vartheta(t) := \int_{[-r,0]} x_{0,\vartheta}(t + u) \, a(\dd u) ,
   \qquad t \in \RR_+ .
 \]
Applying Lemma \ref{res}, we obtain
 \[
   y_\vartheta(t)
   = \sum_{\underset{\SC\Re(\lambda)\geq c}
                    {\lambda \in \Lambda_\vartheta}}
      \int_{[-r,0]}
       \underset{\SC z=\lambda}{\Res}
        \left(\frac{\ee^{z(t+u)}}{h_\vartheta(z)}\right)
       a(\dd u)
     + \int_{[-r,0]} \psi_{\vartheta,c}(t + u) \, a(\dd u) .
 \]
Here we have
 \[
   \int_{[-r,0]} \psi_{\vartheta,c}(t + u) \, a(\dd u) = \oo(\ee^{ct}) \qquad
   \text{as \ $t \to \infty$.}
 \]
Indeed,
 \[
   \lim_{t\to\infty}
    \ee^{-ct} \int_{[-r,0]} \psi_{\vartheta,c}(t + u) \, a(\dd u)
   = \lim_{t\to\infty}
      \int_{[-r,0]}
       [\ee^{-c(t+u)} \psi_{\vartheta,c}(t + u)] \ee^{cu} \, a(\dd u)
   = 0 ,
 \]
 since \ $a$ \ is a finite signed measure on \ $[-r, 0]$.
\ Moreover,
 \begin{align*}
  \underset{\SC z=\lambda}{\Res}
   \left(\frac{\ee^{z(t+u)}}{h_\vartheta(z)}\right)
  &= \ee^{\lambda(t+u)}
     \sum_{k=-m_\vartheta(\lambda)}^{-1}
      \frac{A_{\vartheta,k}(\lambda)}{(-1-k)!} \, (t+u)^{-1-k} \\
  &= \ee^{\lambda(t+u)}
     \sum_{k=-m_\vartheta(\lambda)}^{-1}
      A_{\vartheta,k}(\lambda)
      \sum_{\ell=0}^{-1-k}
       \frac{t^\ell u^{-1-k-\ell}}{\ell!\,(-1-k-\ell)!} \\
  &= \ee^{\lambda(t+u)}
     \sum_{\ell=0}^{m_\vartheta(\lambda)-1}
      \frac{t^\ell}{\ell!}
      \sum_{k=-m_\vartheta(\lambda)}^{-1-\ell}
       \frac{A_{\vartheta,k}(\lambda)}{(-1-k-\ell)!} \, u^{-1-k-\ell} \\
  &= \ee^{\lambda t}
     \sum_{\ell=0}^{m_\vartheta(\lambda)-1}  
      \frac{t^\ell}{\ell!} \,
      \underset{\SC z=\lambda}{\Res}
       \left(\frac{(z-\lambda)^\ell \ee^{zu}}{h_\vartheta(z)}\right) .
 \end{align*}
Consequently, we obtain for each \ $\vartheta \in \RR$ \ and each
 \ $c \in \RR$, \ the representation
 \begin{equation}\label{intres}
  y_\vartheta(t)
  = \sum_{\underset{\SC\Re(\lambda)\geq c}
                   {\lambda \in \Lambda_\vartheta}}
     P_{\vartheta,\lambda}(t) \, \ee^{\lambda t}
    + \Psi_{\vartheta,c}(t) \qquad \text{as \ $t \to \infty$,}
 \end{equation}
 where \ $\Psi_{\vartheta,c} : \RR_+ \to \RR$ \ is a continuous function with
 \ $\Psi_{\vartheta,c}(t) = \oo(\ee^{ct})$ \ as \ $t \to \infty$.
\ Hence we need to analyse the asymptotic behavior of the right hand side of
 \eqref{Yvartheta} as \ $T \to \infty$, \ replacing \ $y_\vartheta(t)$ \ by 
 \ $P_{\vartheta,\lambda}(t) \ee^{\lambda t}$.

First we derive a good estimate for the second term of the right hand side of
 \eqref{Yvartheta}.

\begin{Lem}\label{CS}
Let \ $(y(t))_{t\in\RR_+}$ \ be a continuous deterministic function.
Let \ $a$ \ be a signed measure on \ $[-r, 0]$.
\ Put
 \[
   I(t) := \int_{[-r,0]} \int_u^0 y(t + u - s) X_0(s) \, \dd s \, a(\dd u) ,
   \qquad t \in [r, \infty) .
 \]
Then for each \ $T \in [r, \infty)$,
 \begin{gather}\label{CS1}
  I(T)^2
  \leq |a|([-r,0]) \int_{-r}^0 X_0(s)^2 \, \dd s \int_0^T y(v)^2 \, \dd v
  \leq \|a\| \int_{-r}^0 X_0(s)^2 \, \dd s \int_0^T y(v)^2 \, \dd v , \\
  \begin{aligned}
   \int_r^T I(t)^2 \, \dd t
   &\leq \int_{[-r,0]} (-u) \, |a|(\dd u) \int_{-r}^0 X_0(s)^2 \, \dd s
         \int_0^T y(v)^2 \, \dd v \\
   &\leq r \|a\| \int_{-r}^0 X_0(s)^2 \, \dd s \int_0^T y(v)^2 \, \dd v ,
  \end{aligned} \label{CS2}
 \end{gather}
 where \ $|a|$ \ and \ $\|a\| := |a|([-r,0])$ \ denotes the variation and the
 total variation of the signed measure \ $a$, \ respectively.
\end{Lem}

\noindent{\bf Proof.}
For each \ $t \in [r, \infty)$, \ by Fubini's theorem, 
 \[
   I(t)
   = \int_{-r}^0 X_0(s) \int_{[-r,s]} y(t + u - s) \, a(\dd u) \, \dd s .
 \]
By the Cauchy--Schwarz inequality,
 \begin{align*}
  I(t)^2
  &\leq \int_{-r}^0 X_0(s)^2 \, \dd s
        \int_{-r}^0
         \left(\int_{[-r,s]} y(t + u - s) \, a(\dd u)\right)^2 \dd s \\
  &\leq \int_{-r}^0 X_0(s)^2 \, \dd s
        \int_{-r}^0 \int_{[-r,s]} y(t + u - s)^2 \, |a|(\dd u) \, \dd s ,
 \end{align*}
 where, by Fubini's theorem,  
 \begin{align*}
  \int_{-r}^0 \int_{[-r,s]} y(t + u - s)^2 \, |a|(\dd u) \, \dd s
  &= \int_{[-r,0]} \int_u^0 y(t + u - s)^2 \, \dd s \, |a|(\dd u) \\
  &= \int_{[-r,0]} \int_{t+u}^t y(v)^2 \, \dd v \, |a|(\dd u)
   \leq \int_{[-r,0]} |a|(\dd u) \int_0^t y(v)^2 \, \dd v ,
 \end{align*}
 hence we obtain \eqref{CS1}.
Moreover,
 \[
   \int_r^T I(t)^2 \, \dd t
   \leq \int_{-r}^0 X_0(s)^2 \, \dd s
        \int_r^T \int_{-r}^0 \int_{[-r,s]}
         y(t + u - s)^2 \, |a|(\dd u) \, \dd s \, \dd t ,   
 \]
 where, by Fubini's theorem,  
 \begin{align*}
  &\int_r^T \int_{-r}^0 \int_{[-r,s]}
    y(t + u - s)^2 \, |a|(\dd u) \, \dd s \, \dd t
   = \int_{-r}^0 \int_{[-r,s]} \int_r^T
      y(t + u - s)^2 \, \dd t \, |a|(\dd u) \, \dd s \\
  &\qquad\qquad
   = \int_{-r}^0 \int_{[-r,s]} \int_{r+u-s}^{T+u-s} y(v)^2 \, \dd v \, |a|(\dd u) \, \dd s
   \leq \int_0^T y(v)^2 \, \dd v \int_{-r}^0 \int_{[-r,s]} \, |a|(\dd u) \, \dd s \\
  &\qquad\qquad
   = \int_0^T y(v)^2 \, \dd v \int_{[-r,0]} \int_u^0 \, \dd s \, |a|(\dd u)
   = \int_0^T y(v)^2 \, \dd v \int_{[-r,0]} (-u) \, |a|(\dd u) ,
 \end{align*}
 hence we obtain \eqref{CS2}.
\proofend
 
\begin{Lem}\label{erg}
Let \ $(y(t))_{t\in\RR_+}$ \ be a continuous deterministic
 function with \ $\int_0^\infty y(t)^2 \, \dd t < \infty$.
\ Let \ $\vartheta \in \RR$.  
\ Suppose that \ $(Y(t))_{t\in\RR_+}$ \ is a continuous stochastic process such
 that
 \begin{equation}\label{Y}
  Y(t)
  = y(t) X_0(0)
    + \vartheta \int_{[-r,0]} \int_v^0
       y(t + v - s) X_0(s) \, \dd s \, a(\dd v)
    + \int_0^t y(t - s) \, \dd W(s) 
 \end{equation}
 for \ $t \in [r,\infty)$.
\ Then
 \begin{gather}
  \frac{1}{T} \int_0^T Y(t) \, \dd t \stoch 0 \qquad
   \text{as \ $T \to \infty$,} \label{erg1} \\
  \frac{1}{T} \int_0^T Y(t)^2 \, \dd t
  \stoch \int_0^\infty y(t)^2 \, \dd t \qquad
  \text{as \ $T \to \infty$.} \label{erg2}
 \end{gather}
\end{Lem}

\noindent{\bf Proof.}
Applying Lemma 4.3 in Gushchin and K\"uchler \cite{GusKuc1999} for the special
 case \ $X_0(s) = 0$, \ $s \in [-1, 0]$, \ we obtain
 \begin{gather*}
  \frac{1}{T} \int_0^T \int_0^t y(t - s) \, \dd W(s) \, \dd t \stoch 0 \qquad
  \text{as \ $T \to \infty$,} \\
  \frac{1}{T} \int_0^T \left(\int_0^t y(t - s) \, \dd W(s)\right)^2 \dd t
  \stoch \int_0^\infty y(t)^2 \, \dd t \qquad
  \text{as \ $T \to \infty$.}
 \end{gather*}
We have
 \[
   \frac{1}{T} \int_0^T Y(t) \, \dd t
   = \frac{1}{T} \int_0^r Y(t) \, \dd t
     + X_0(0) I_1(T)
     + \frac{\vartheta}{T} \int_r^T Z(t) \, \dd t
     + \frac{1}{T} \int_r^T \int_0^t y(t - s) \, \dd W(s) \, \dd t
 \]
 for \ $T \in \RR_+$ \ with
 \begin{gather}\nonumber
  I_1(T) := \frac{1}{T} \int_r^T y(t) \, \dd t , \qquad T \in \RR_+ , \\
  Z(t) := \int_{[-r,0]} \int_u^0 y(t + u - s) X_0(s) \, \dd s \, a(\dd u) ,
  \qquad t \in [r, \infty) . \label{Z}
 \end{gather}
By the Cauchy--Schwarz inequality and by \eqref{CS2},
 \begin{gather*}
  |I_1(T)| 
  \leq \sqrt{\int_r^T \frac{1}{T^2} \, \dd t
             \int_r^T y(t)^2 \, \dd t}
  \leq \sqrt{\frac{1}{T} \int_0^\infty y(t)^2 \, \dd t}
  \to 0 , \\
  \begin{aligned}
   \left|\frac{1}{T} \int_r^T Z(t) \, \dd t\right|
   &\leq \sqrt{\frac{1}{T} \int_r^T Z(t)^2 \, \dd t}
    \leq \sqrt{\frac{r\|a\|}{T} \int_{-r}^0 X_0(s)^2 \, \dd s
               \int_0^\infty y(v)^2 \, \dd v}
    \as 0
  \end{aligned}
 \end{gather*}
 as \ $T \to \infty$, \ hence we obtain \eqref{erg1}.
Moreover,
 \[
   \frac{1}{T} \int_0^T Y(t)^2 \, \dd t
   = \frac{1}{T} \int_0^r Y(t)^2 \, \dd t
     + I_2(T) + 2 I_3(T)
     + \frac{1}{T} \int_r^T \left(\int_0^t y(t - s) \, \dd W(s)\right)^2 \dd t
 \]
 for \ $T \in \RR_+$, \ with
 \begin{gather*}
  I_2(T) := \frac{1}{T} \int_r^T (y(t) X_0(0) + \vartheta Z(t))^2 \, \dd t , \\
  I_3(T) := \frac{1}{T}
                \int_r^T 
                 (y_i(t) X_0(0) + \vartheta Z_i(t))
                 \left(\int_0^t y_i(t - s) \, \dd W(s)\right) \dd t .
 \end{gather*}
Again by \eqref{CS2},
 \begin{align*}
  0 &\leq I_2(T)
     \leq \frac{1}{T}
        \int_r^T
         2 (y(t)^2 X_0(0)^2 + \vartheta^2 Z(t)^2) \, \dd t \\
  &\leq \frac{2X_0(0)^2}{T} \int_0^\infty y(t)^2 \, \dd t
        + \frac{2}{T}
          \vartheta^2 r \|a\| \int_{-r}^0 X_0(s)^2 \, \dd s
          \int_0^\infty y(v)^2 \, \dd v
  \as 0
 \end{align*}
 as \ $T \to \infty$, \ and
 \begin{align*}
  |I_3(T)| 
  &\leq \frac{2}{T} 
        \sqrt{\int_r^T (y(t) X_0(0) + \vartheta Z(t))^2 \, \dd t 
              \int_r^T \left(\int_0^t y(t - s) \, \dd W(s)\right)^2 \dd t} \\
  &= 2 \sqrt{\frac{I_2(T)}{T} 
     \int_r^T \left(\int_0^t y(t - s) \, \dd W(s)\right)^2 \dd t}
  \stoch 0 \qquad \text{as \ $T \to \infty$,}
 \end{align*}
 hence we obtain \eqref{erg2}.
\proofend

\begin{Lem}\label{super}
Let \ $y_\ell(t) = t^{\alpha_\ell} \Re(c_\ell \, \ee^{\lambda_\ell t})$,
 \ $t \in \RR_+$, \ $\ell \in \{1, 2\}$, \ with some \ $\alpha_\ell \in \ZZ_+$,
 \ $c_\ell, \lambda_\ell \in \CC$ \ with \ $\Re(\lambda_\ell) \in \RR_{++}$,
 \ $\ell \in \{1, 2\}$. 
\ Let \ $\vartheta \in \RR$.  
\ Suppose that \ $(Y_\ell(t))_{t\in\RR_+}$, \ $\ell \in \{1, 2\}$, \ are
 continuous stochastic processes admitting representation \eqref{Y} on
 \ $[r, \infty)$ \ with \ $y = y_\ell$, \ $\ell \in \{1, 2\}$, \ respectively.
Then
 \begin{equation}\label{super1}
  t^{-\alpha_1} \ee^{-t\Re(\lambda_1)} Y_1(t)
  - \Re\bigl(c_1 U^{(\vartheta)}_{\lambda_1}
              \ee^{\ii t\Im(\lambda_1)}\bigr)
  \as 0 , \qquad \text{as \ $t \to \infty$,} 
 \end{equation}
 and
 \begin{equation}\label{super2}
  \begin{aligned}
   &T^{-\alpha_1-\alpha_2} \ee^{-T\Re(\lambda_1+\lambda_2)}
    \int_0^T Y_1(t) Y_2(t) \, \dd t \\
   &- \int_0^\infty
       \ee^{-t\Re(\lambda_1+\lambda_2)}
       \Re\bigl(c_1 U^{(\vartheta)}_{\lambda_1}
                \ee^{\ii(T-t)\Im(\lambda_1)}\bigr)
       \Re\bigl(c_2 U^{(\vartheta)}_{\lambda_2}
                \ee^{\ii(T-t)\Im(\lambda_2)}\bigr) \,
       \dd t
    \as 0 ,
  \end{aligned}
 \end{equation}
 as \ $T \to \infty$.
\ Particularly, if \ $c_1, c_2 \in \RR$ \ and
 \ $\lambda_1, \lambda_2 \in \RR_{++}$, \ then
 \begin{gather*}
  t^{-\alpha_1} \ee^{-\lambda_1 t} Y_1(t)
  \as c_1 U^{(\vartheta)}_{\lambda_1} ,
  \qquad \text{as \ $t \to \infty$,} \\
  T^{-\alpha_1-\alpha_2} \ee^{-T(\lambda_1+\lambda_2)}
  \int_0^T Y_1(t) Y_2(t) \, \dd t
  \as \frac{c_1 c_2}{\lambda_1+\lambda_2} \,
      U^{(\vartheta)}_{\lambda_1} U^{(\vartheta)}_{\lambda_1} ,
  \qquad \text{as \ $T \to \infty$.}
 \end{gather*}
\end{Lem}

\noindent{\bf Proof.}
Note that
 \[
   t^{-\alpha_1} \ee^{-t\Re(\lambda_1)} Y_1(t)
   - \Re\bigl(c_1 U^{(\vartheta)}_{\lambda_1} \ee^{\ii t\Im(\lambda_1)}\bigr)
   = - I_1(t) + I_2(t) - I_3(t) , \qquad t \in [r, \infty) ,
 \]
 with
 \begin{align*}
  I_1(t)
  &:= \vartheta
      \int_{[-r,0]} \int_v^0
       \left[1 - \left(1 - \frac{s-v}{t}\right)^{\alpha_1}\right]
       \Re\bigl(c_1 \ee^{\ii t\Im(\lambda_1)-\lambda_1(s-v)}\bigr) X_0(s)
       \, \dd s \, a(\dd v) , \\
  I_2(t)
  &:= \int_0^t
       \left[\left(1 - \frac{s}{t}\right)^{\alpha_1} - 1\right]
       \Re\bigl(c_1 \ee^{\ii t\Im(\lambda_1)-\lambda_1 s}\bigr)
       \, \dd W(s) , \\
  I_3(t)
  &:= \int_t^\infty 
       \Re\bigl(c_1 \ee^{\ii t\Im(\lambda_1)-\lambda_1 s}\bigr) \, \dd W(s) .
 \end{align*}
By the dominated convergence theorem, \ $I_1(t) \as 0$ \ as \ $t \to \infty$. 
\ Moreover,
 \[
   I_2(t)
   = \sum_{k=1}^{\alpha_1}
      (-1)^k \binom{\alpha_1}{k} t^{-k}
      \int_0^t
       s^k \Re\bigl(c_1 \ee^{\ii t\Im(\lambda_1)-\lambda_1 s}\bigr)
       \, \dd W(s)
   \as 0 \qquad \text{as \ $t \to \infty$}
 \]
 by the strong law of martingales, see, e.g., Liptser and Shiryaev
 \cite[Chapter 2, \S 6, Theorem 10]{LipShiI}.
Obviously, \ $I_3(t) \as 0$ \ as \ $t \to \infty$, \ hence we obtain
 \eqref{super1}.
 
In order to prove \eqref{super2}, put
 \[
   V_\ell(t) := \Re\bigl(c_\ell \, U^{(\vartheta)}_{\lambda_\ell}
                         \ee^{\ii t\Im(\lambda_\ell)}\bigr) , \qquad
   t \in \RR , \quad \ell \in \{1, 2\} .
 \]
For each \ $T \in \RR_+$, \ we have
 \[
   \int_0^T \ee^{-t\Re(\lambda_1+\lambda_2)} V_1(T-t) V_2(T-t) \, \dd t
   = \ee^{-T\Re(\lambda_1+\lambda_2)}
     \int_0^T \ee^{t\Re(\lambda_1+\lambda_2)} V_1(t) V_2(t) \, \dd t ,
 \]
 hence
 \begin{align*}
  &T^{-\alpha_1-\alpha_2} \ee^{-T\Re(\lambda_1+\lambda_2)}
   \int_0^T Y_1(t) Y_2(t) \, \dd t
   - \int_0^\infty
      \ee^{-t\Re(\lambda_1+\lambda_2)} V_1(T-t) V_2(T-t) \, \dd t \\
  &= J_1(T) + J_2(T) + J_3(t) - J_4(T) - J_5(T) 
 \end{align*}
 with
 \begin{align*}
  J_1(T)
  &:= T^{-\alpha_1-\alpha_2} \ee^{-T\Re(\lambda_1+\lambda_2)}
      \int_0^T
       [Y_1(t) - t^{\alpha_1} \ee^{t\Re(\lambda_1)} V_1(t)]
       [Y_2(t) - t^{\alpha_2} \ee^{t\Re(\lambda_2)} V_2(t)]
       \, \dd t , \\
  J_2(T)
  &:= T^{-\alpha_1-\alpha_2} \ee^{-T\Re(\lambda_1+\lambda_2)}
      \int_0^T
       [Y_1(t) - t^{\alpha_1} \ee^{t\Re(\lambda_1)} V_1(t)]
       t^{\alpha_2} \ee^{t\Re(\lambda_2)} V_2(t) \, \dd t, \\
  J_3(T)
  &:= T^{-\alpha_1-\alpha_2} \ee^{-T\Re(\lambda_1+\lambda_2)}
      \int_0^T
       [Y_2(t) - t^{\alpha_2} \ee^{t\Re(\lambda_2)} V_2(t)]
       t^{\alpha_1} \ee^{t\Re(\lambda_1)} V_1(t) \, \dd t, \\
  J_4(T)
  &:= \ee^{-T\Re(\lambda_1+\lambda_2)}
      \int_0^T
       \biggl(1 - \frac{t^{\alpha_1+\alpha_2}}{T^{\alpha_1+\alpha_2}}\biggr)
       \ee^{t\Re(\lambda_1+\lambda_2)} V_1(t) V_2(t) \, \dd t , \\
  J_5(T)
  &:= \int_T^\infty
       \ee^{-t\Re(\lambda_1+\lambda_2)} V_1(T-t) V_2(T-t) \, \dd t .
 \end{align*}
By \eqref{super1} and L'H\^ospital's rule, \ $J_1(T) \as 0$ \ as
 \ $T \to \infty$.
\ By the Cauchy--Schwarz inequality, \ $|J_2(T)| \leq \sqrt{J_6(T) J_7(T)}$,
 \ $T \in \RR_+$, \ where, by \eqref{super1} and L'H\^ospital's rule,
 \[
   J_6(T)
   := T^{-2\alpha_1} \ee^{-2T\Re(\lambda_1)}
      \int_0^T
       [Y_1(t) - t^{\alpha_1} \ee^{t\Re(\lambda_1)} V_1(t)]^2
       \, \dd t
   \as 0 , \qquad \text{as \ $T \to \infty$,}
 \]
 and
 \begin{align*}
  J_7(T)
  &:= T^{-2\alpha_2} \ee^{-2T\Re(\lambda_2)}
      \int_0^T t^{2\alpha_2} \ee^{2t\Re(\lambda_2)} V_2(t)^2 \, \dd t \\
  &= \int_0^T
      \biggl(1 - \frac{t}{T}\biggr)^{2\alpha_2}
      \ee^{-2t\Re(\lambda_2)} V_2(T - t)^2 \, \dd t
   \leq \frac{1}{2\Re(\lambda_2)}
        \sup_{t\in\RR} V_2(t)^2
   < \infty \quad \text{a.s.,}
 \end{align*}
 since \ $(V_2(t))_{t\in\RR}$ \ is a continuous and periodic process.
Consequently, \ $J_2(T) \as 0$ \ as \ $T \to \infty$.
\ In a similar way, \ $J_3(T) \as 0$ \ as \ $T \to \infty$, \ and
 \begin{align*}
  J_4(T)
  &\leq \left(\sup_{t\in\RR} |V_1(t) V_2(t)|\right)
        \ee^{-T\Re(\lambda_1+\lambda_2)}
        \int_0^T
         \left(1 - \frac{t^{\alpha_1+\alpha_2}}{T^{\alpha_1+\alpha_2}}\right)
         \ee^{t\Re(\lambda_1+\lambda_2)} \, \dd t \\
  &= \biggl(\sup_{t\in\RR} |V_1(t) V_2(t)|\biggr)
     \int_0^T
      \biggl[1 - \biggl(1 - \frac{t}{T}\biggr)^{\alpha_1+\alpha_2}\biggr]
      \ee^{-t\Re(\lambda_1+\lambda_2)} \, \dd t
   \as 0
 \end{align*}
 as \ $T \to \infty$ \ by the dominated convergence theorem.
Finally,
 \[
   |J_5(T)|
   \leq \frac{\ee^{-T\Re(\lambda_1+\lambda_2)}}{\Re(\lambda_1+\lambda_2)}
        \sup_{t\in\RR} |V_1(t) V_2(t)|
   \as 0 \qquad \text{as \ $T \to \infty$,}
 \]
 hence we obtain \eqref{super2}.  
\proofend

\noindent{\bf Proof of Theorem \ref{LAN}.}
The continuous process \ $\bigl(Y^{(\vartheta)}(t)\bigr)_{t\in\RR_+}$
 \ admits the representation \eqref{Yvartheta} on \ $[r,\infty)$.
\ The aim of the following discussion is to show that the function
 \ $(y_\vartheta(t))_{t\in\RR_+}$ \ is square integrable. 
Let \ $c \in (v_\vartheta^*,  0)$, \ and apply the representation \eqref{intres}.
By the definition of \ $v_\vartheta^*$, \ we obtain \ $P_{\vartheta,\lambda} = 0$
 \ for each \ $\lambda \in \Lambda_\vartheta$ \ with
 \ $\Re(\lambda) > v_\vartheta^*$, \ and hence for each
 \ $\lambda \in \Lambda_\vartheta$ \ with \ $\Re(\lambda) \geq c$.
\ Thus the representation \eqref{intres} gives
 \ $y_\vartheta(t) = \oo(\ee^{ct})$ \ as \ $t \to \infty$. 
\ Since \ $(y_\vartheta(t))_{t\in\RR_+}$ \ is continuous, the function
 \ $(\ee^{-ct} y_\vartheta(t))_{t\in\RR_+}$ \ is bounded, implying
 \ $\int_0^\infty y_\vartheta(t)^2 \, \dd t < \infty$.       
\ Hence we can apply Lemma \ref{erg} to obtain
 \[
   J_{\vartheta,T}
   = \frac{1}{T} \int_0^T Y^{(\vartheta)}(t)^2 \, \dd t
   = \frac{1}{T} \int_0^r Y^{(\vartheta)}(t)^2 \, \dd t
     + \frac{1}{T} \int_r^T Y^{(\vartheta)}(t)^2 \, \dd t
   \stoch
   \int_0^\infty y_\vartheta(t)^2 \, \dd t
   = J_\vartheta
 \]
 as \ $T \to \infty$, \ where \ $J_\vartheta > 0$, \ since
 \ $x_{0,\vartheta} \ne 0$ \ and \ $a \ne 0$ \ implies \ $y_\vartheta \ne 0$.

In case of \ $\vartheta = 0$, \ we have \ $h_0(\lambda) = \lambda$,
 \ $\Lambda_0 = \{0\}$, \ $m_0(0) = 1$ \ and
 \ $P_{0,0}(t) = A_{0,-1}(0) \int_{[-r,0]} a(\dd u) = a([-r,0])$,
 \ $t \in \RR$, \ since \ $1/h_0(z) = z^{-1}$ \ yields \ $A_{0,-1}(0) = 1$.
\ The assumption \ $a([-r,0]) = 0$ \ implies \ $P_{0,0} = 0$, \ and hence
 \ $v_0^* = -\infty$ \ and \ $m_0^* = -\infty$. 
\ Moreover, the assumption \ $a([-r,0]) = 0$ \ yields
 \[
   y_0(t)
   = \int_{[-r,0]} x_{0,0}(t + u) \, a(\dd u)
   = \begin{cases}
      0 & \text{if \ $t \in [r, \infty)$,} \\
      a([-t,0]) & \text{if \ $t \in [0, r]$,}
     \end{cases}
 \]
 and we obtain the formula for \ $J_0$.

Further, the process 
 \[
   M^{(\vartheta)}(T) := \int_0^T Y^{(\vartheta)}(t) \, \dd W(t) ,
   \qquad T \in \RR_+ ,
 \]
 is a continuous martingale with \ $M^{(\vartheta)}(0) = 0$ \ and with
 quadratic variation
 \[
   \langle M^{(\vartheta)} \rangle(T)
   = \int_0^T Y^{(\vartheta)}(t)^2 \, \dd t ,
 \]
 hence Theorem VIII.5.42 of Jacod and Shiryaev \cite{JSh} yields the statement.
\proofend
 
\noindent{\bf Proof of Theorem \ref{LAQ}.}
For each \ $T \in \RR_{++}$, \ we have
 \begin{gather*}
  \Delta_{\vartheta,T}
  = \frac{1}{T^{m_\vartheta^*+1}} 
    \int_0^T Y^{(\vartheta)}(t) \, \dd W(t) , \qquad
  J_{\vartheta,T}
  = \frac{1}{T^{2(m_\vartheta^*+1})}
    \int_0^T Y^{(\vartheta)}(t)^2 \, \dd t .
 \end{gather*}
The continuous process \ $\bigl(Y^{(\vartheta)}(t)\bigr)_{t\in\RR_+}$
 \ admits the representation \eqref{Yvartheta} on \ $[r,\infty)$.
\ We choose \ $c < 0$ \ with
 \ $c > \sup\{\Re(\lambda)
              : \lambda \in \Lambda_\vartheta , \; \Re(\lambda) < 0\}$,
 \ and apply the representation \eqref{intres}.
The assumption \ $v_\vartheta^* = 0$ \ yields that \ $P_{\vartheta,\lambda} = 0$
 \ for each \ $\lambda \in \Lambda_\vartheta$ \ with \ $\Re(\lambda) > 0$,
 \ hence we obtain
 \begin{equation}\label{yres}
  y_\vartheta(t)
  = \sum_{\lambda \in \Lambda_\vartheta \cap (\ii\RR)}
     P_{\vartheta,\lambda}(t) \, \ee^{\ii t\Im(\lambda)}
    + \Psi_{\vartheta,c}(t) , \qquad t \in \RR_+ .
 \end{equation}
The leading term of the polynomial \ $P_{\vartheta,\lambda}$ \ is
 \ $c_{\vartheta,\lambda,m_\vartheta^*} t^{\tm_\vartheta(\lambda)}$, \ thus, by
 the representation \eqref{Yvartheta},
 \begin{equation}\label{tYres}
  Y^{(\vartheta)}(t)
  = \sum_{\underset{\SC\tm_\vartheta(\lambda)=m_\vartheta^*}
                   {\lambda \in \Lambda_\vartheta \cap (\ii\RR)}}
     c_{\vartheta,\lambda,m_\vartheta^*} \, \ee^{\ii t\Im(\lambda)}
     \int_0^t (t - s)^{m_\vartheta^*} \ee^{-\ii s\Im(\lambda)} \, \dd W(s)
    + \tY(t) , \qquad t \in \RR_+ ,
 \end{equation}
 where \ $(\tY(t))_{t\in\RR_+}$ \ is a continuous process satisfying
 \[
   \tY(t)
   := y_\vartheta(t) X_0(0) + \vartheta I_\vartheta(t)
      + I_\Psi(t) + S(t) ,
   \qquad t \in [r, \infty) ,
 \]
 with
 \begin{gather*}
  I_\vartheta(t) := \int_{[-r,0]} \int_v^0
                     y_\vartheta(t + v - s) X_0(s)
                     \, \dd s \, a(\dd v) , \qquad
  I_\Psi(t) := \int_0^t \Psi_{\vartheta,c}(t - s) \, \dd W(s) , \\
  S(t) := \sum_{\lambda \in \Lambda_\vartheta \cap (\ii\RR)}
           \sum_{\ell=0}^{\tm_\vartheta(\lambda)\land(m_\vartheta^*-1)}
            c_{\vartheta,\lambda,\ell} I_\ell(t) , \qquad
  I_\ell(t) := \int_0^t (t - s)^\ell \ee^{\lambda(t-s)} \, \dd W(s) ,        
 \end{gather*}
 for \ $t \in [r, \infty)$ \ and \ $\ell \in \ZZ_+$.
\ The aim of the following discussion is to show that 
 \ $T^{-2(m_\vartheta^*+1)} \int_r^T \tY(t)^2 \, \dd t \stoch 0$ \ as
 \ $T \to \infty$. 
\ First we show that
 \ $T^{-2(m_\vartheta^*+1)} \int_r^T y_\vartheta(t)^2 \, \dd t \to 0$
 \ as \ $T \to \infty$. 
\ For each \ $\lambda \in \Lambda_\vartheta$ \ with \ $\Re(\lambda) = 0$, \ we
 have
 \[
   T^{-2(m_\vartheta^*+1)} 
   \int_r^T |P_{\vartheta,\lambda}(t) \ee^{\lambda t}|^2 \, \dd t
   = T^{-2(m_\vartheta^*+1)}
     \int_r^T |P_{\vartheta,\lambda}(t)|^2 \, \dd t
   \to 0 \qquad \text{as \ $T \to \infty$,}
 \]
 since the degree of the complex-valued polynomial
 \ $P_{\vartheta,\lambda}$ \ is \ $\tm_\vartheta(\lambda)$, \ which is at most
 \ $m_\vartheta^*$.
\ Moreover,
 \[
   T^{-2(m_\vartheta^*+1)} \int_r^T \Psi_{\vartheta,c}(t)^2 \, \dd t
   \leq T^{-2} \int_r^T \Psi_{\vartheta,c}(t)^2 \, \dd t
   \leq T^{-2} \int_r^\infty \Psi_{\vartheta,c}(t)^2 \, \dd t
   \to 0
 \]
 as \ $T \to \infty$, \ since
 \ $\int_r^\infty |\Psi_{\vartheta,c}(t)|^2 \, \dd t < \infty$, \ because the
 function \ $(\ee^{-ct} |\Psi_{\vartheta,c}(t)|)_{t\in\RR_+}$ \ is bounded.
Thus, by the representation \eqref{yres}, we obtain
 \ $T^{-2(m_\vartheta^*+1)} \int_r^T y_\vartheta(t)^2 \, \dd t \to 0$ \ as
 \ $T \to \infty$. 
\ Next, by Lemma \ref{CS},
 \[
   \int_r^T I_\vartheta(t)^2 \, \dd t
   \leq r \|a\| \int_{-r}^0 X_0(s)^2 \, \dd s
        \int_0^T y_\vartheta(v)^2 \, \dd v ,
 \]
 hence we obtain
 \ $T^{-2(m_\vartheta^*+1)} \int_r^T I_\vartheta(t)^2 \, \dd t \stoch 0$ \ as
 \ $T \to \infty$. 
\ Further, we have
 \begin{equation}\label{psi}
  \begin{aligned}
   &\EE\left(\int_r^T I_\Psi(t)^2 \, \dd t\right)
   = \int_r^T \left(\int_0^t \Psi_{\vartheta,c}(t - s)^2 \, \dd s\right) \dd t
   = \int_r^T \left(\int_0^t \Psi_{\vartheta,c}(u)^2 \, \dd u\right) \dd t \\
   &= \int_r^T \Psi_{\vartheta,c}(u)^2 \left(\int_u^T \dd t\right) \dd u
   = \int_r^T (T - u) \Psi_{\vartheta,c}(u)^2 \, \dd u
  \leq T \int_r^T \Psi_{\vartheta,c}(u)^2 \, \dd u ,
  \end{aligned}
 \end{equation}
 hence \ $T^{-2} \int_r^T I_\Psi(t)^2 \, \dd t \stoch 0$ \ as
 \ $T \to \infty$. 
\ Finally, for each \ $\lambda \in \Lambda_\vartheta$ \ with
 \ $\Re(\lambda) = 0$ \ and each \ $\ell \in \ZZ_+$ \ with 
 \ $\ell \leq m_\vartheta^* - 1$, \ we have
 \ $T^{-2(m_\vartheta^*+1)} \int_0^T I_\ell(t)^2 \, \dd t \stoch 0$ \ as
 \ $T \to \infty$.
\ Indeed,
 \begin{align*}
  &\frac{1}{T^{2(m_\vartheta^*+1)}} \EE\left(\int_0^T I_\ell(t)^2 \, \dd t\right)
  = \frac{1}{T^{2(m_\vartheta^*+1)}}
    \int_0^T \left(\int_0^t (t - s)^{2\ell} \, \dd s \right) \dd t \\
  &= \frac{1}{T^{2(m_\vartheta^*+1)}}
    \int_0^T \left(\int_0^t u^{2\ell} \, \dd u \right) \dd t
  = \frac{1}{(2\ell+1)(2\ell+2)T^{2(m_\vartheta^*-\ell)}}
  \to 0 \qquad \text{as \ $T \to \infty$.}
 \end{align*}
Hence we conclude
 \ $T^{-2(m_\vartheta^*+1)} \int_0^T \tY(t)^2 \, \dd t \stoch 0$ \ as
 \ $T \to \infty$. 
 
Introduce the complex-valued processes \ $(Z_{\varphi,\ell}(t))_{t\in\RR_+}$,
 \ $\varphi \in \RR$, \ $\ell \in \ZZ_+$, \ by
 \[
   Z_{\varphi,\ell}(t)
   := \int_0^t (t - s)^\ell \ee^{-\ii\varphi s} \, \dd W(s) .
 \]
Note that \ $Z_{0,0} = W$.
\ For each \ $T \in \RR_{++}$, \ consider the complex-valued processes
 \ $(Z^T_{\varphi,\ell}(t))_{s\in[0,1]}$, \ $\varphi \in \RR$,
 \ $\ell \in \ZZ_+$, \ and the real-valued process \ $(X^T(s))_{s\in[0,1]}$,
 \ defined by
 \begin{gather*}
  Z^T_{\varphi,\ell}(s)
  := \frac{1}{T^{\ell+\frac{1}{2}}} Z_{\varphi,\ell}(Ts)
  = \int_0^s (s - u)^\ell \ee^{-\ii T\varphi u} \, \dd W^T(u) , \\
  X^T(s)
  := \sum_{\underset{\SC\tm_\vartheta(\lambda)=m_\vartheta^*}
                    {\lambda \in \Lambda_\vartheta \cap (\ii\RR)}}
      c_{\vartheta,\lambda,m_\vartheta^*}
      \ee^{\ii Ts\Im(\lambda)} Z^T_{\Im(\lambda),m_\vartheta^*}(s) .
 \end{gather*}
Then, for each \ $T \in \RR_{++}$, \ we have
 \[
   Y^{(\vartheta)}(t)
   = T^{m_\vartheta^*+\frac{1}{2}} X^T\left(\frac{t}{T}\right) + \tY(t) , \qquad
   t \in \RR_+ ,
 \]
 and hence,
 \begin{gather*}
  \Delta_{\vartheta,T}
  = \frac{1}{\sqrt{T}} \int_0^T X^T\Bigl(\frac{t}{T}\Bigr) \, \dd W(t)
    + I_1(T)
  = \int_0^1 X^T(s) \, \dd W^T(s) + I_1(T) , \\
  J_{\vartheta,T}
  = \frac{1}{T} \int_0^T \Bigl|X^T\Bigl(\frac{t}{T}\Bigr)\Bigr|^2 \dd t
    + 2 I_2(T) + I_3(T)
  = \int_0^1 |X^T(s)|^2 \, \dd s + 2 I_2(T) + I_3(T) ,
 \end{gather*}
 with
 \begin{gather*}
  I_1(T) := \frac{1}{T^{m_\vartheta^*+1}} \int_0^T \tY(t) \, \dd W(t) , \qquad
  I_2(T) := \frac{1}{T^{m_\vartheta^*+\frac{3}{2}}}
            \int_0^T
             \Re\Bigl(X^T\Bigl(\frac{t}{T}\Bigr) \overline{\tY(t)}\Bigr)
             \, \dd t , \\
  I_3(T) := \frac{1}{T^{2(m_\vartheta^*+1)}} \int_0^T |\tY(t)|^2 \, \dd t .
 \end{gather*}
Introducing the process
 \[
   Y^T(t) := \int_0^t X^T(s) \, \dd W^T(s) , 
   \qquad t \in \RR_+ , \qquad T \in \RR_{++} ,
 \]
 we have
 \[
   \int_0^t X^T(s)^2 \, \dd s = [Y^T, Y^T](t) ,
   \qquad t \in \RR_+ , \qquad T \in \RR_{++} ,
 \]
 where \ $([U, V](t))_{t\in\RR_+}$ \ denotes the quadratic covariation process
 of the processes \ $(U(t))_{t\in\RR_+}$ \ and \ $(V(t))_{t\in\RR_+}$.
\ Moreover, 
 \begin{align*}
  Y^T(t)
  &= \sum_{\underset{\SC\tm_\vartheta(\lambda)=m_\vartheta^*}
                    {\lambda \in \Lambda_\vartheta \cap (\ii\RR)}}
      c_{\vartheta,\lambda,m_\vartheta^*}
      \int_0^t
       \ee^{\ii Ts\Im(\lambda)} Z^T_{\Im(\lambda),m_\vartheta^*}(s)
       \, \dd W^T(s) \\
  &= \sum_{\underset{\SC\tm_\vartheta(\lambda)=m_\vartheta^*}
                    {\lambda \in \Lambda_\vartheta \cap (\ii\RR)}}
      c_{\vartheta,\lambda,m_\vartheta^*}
      \int_0^t
       Z^T_{\Im(\lambda),m_\vartheta^*}(s)
       \, \dd \overline{Z^T_{\Im(\lambda),0}(s)} .
 \end{align*}
By the functional central limit theorem,
 \[
   \bigl(Z^T_{\Im(\lambda),0} 
         : \lambda \in \Lambda_\vartheta\cap(\ii\RR)\bigr)
   \distr 
   \bigl(\cZ_{\Im(\lambda),0}
         : \lambda \in \Lambda_\vartheta\cap(\ii\RR)\bigr) , \qquad
   \text{as \ $T \to \infty$,}
 \]
 since
 \[
   Z^T_{\varphi,0}(s)
   = \int_0^s \cos(T\varphi u) \, \dd W(u)
     - \ii \int_0^s \sin(T\varphi u) \, \dd W(u) , \qquad
   Z^T_{-\varphi,0}(s) = \overline{Z^T_{\varphi,0}(s)} ,
 \]
 for all \ $s \in [0, 1]$, \ $\varphi \in \RR$ \ and \ $T \in \RR_{++}$, \ and
 \begin{gather*}
  \begin{aligned}
   \int_0^s \cos(T\varphi_1 u) \cos(T\varphi_2 u) \, \dd u
   &= \frac{1}{2}
      \int_0^s
       [\cos(T(\varphi_1-\varphi_2)u) + \cos(T(\varphi_1+\varphi_2)u)]
       \, \dd u \\
   &= \begin{cases}
       \frac{s}{2} + \frac{\sin(2T\varphi_1 s)}{4T\varphi_1}
       \to \frac{s}{2} ,
        & \text{if \ $\varphi_1 = \varphi_2$,} \\[1mm]
       \frac{\sin(T(\varphi_1-\varphi_2)s)}{2T(\varphi_1-\varphi_2)}
       + \frac{\sin(T(\varphi_1+\varphi_2)s)}{2T(\varphi_1+\varphi_2)}
       \to 0 ,
        & \text{if \ $\varphi_1 \ne \varphi_2$,}
      \end{cases}
  \end{aligned} \\
  \begin{aligned}
   \int_0^s \sin(T\varphi_1 u) \sin(T\varphi_2 u) \, \dd u
   &= \frac{1}{2}
      \int_0^s
       [\cos(T(\varphi_1-\varphi_2)u) - \cos(T(\varphi_1+\varphi_2)u)]
       \, \dd u \\
   &= \begin{cases}
       \frac{s}{2} - \frac{\sin(2T\varphi_1 s)}{4T\varphi_1}
       \to \frac{s}{2} ,
        & \text{if \ $\varphi_1 = \varphi_2$,} \\[1mm]
       \frac{\sin(T(\varphi_1-\varphi_2)s)}{2T(\varphi_1-\varphi_2)}
       - \frac{\sin(T(\varphi_1+\varphi_2)s)}{2T(\varphi_1+\varphi_2)}
       \to 0 ,
        & \text{if \ $\varphi_1 \ne \varphi_2$,}
      \end{cases}
  \end{aligned} \\
  \begin{aligned}
   \int_0^s \sin(T\varphi_1 u) \cos(T\varphi_2 u) \, \dd u
   &= \frac{1}{2}
      \int_0^s
       [\sin(T(\varphi_1+\varphi_2)u) + \sin(T(\varphi_1-\varphi_2)u)]
       \, \dd u \\
   &= \begin{cases}
       \frac{1-\cos(2T\varphi_1 s)}{4T\varphi_1}
       \to 0 ,
        & \text{if \ $\varphi_1 = \varphi_2$,} \\[1mm]
       \frac{1-\cos(T(\varphi_1+\varphi_2)s)}{2T(\varphi_1+\varphi_2)}
       + \frac{1-\cos(T(\varphi_1-\varphi_2)s)}{2T(\varphi_1-\varphi_2)}
       \to 0 ,
        & \text{if \ $\varphi_1 \ne \varphi_2$,}
      \end{cases}
  \end{aligned}
 \end{gather*}
 as \ $T \to \infty$ \ for all \ $s \in [0, 1]$ \ and
 \ $\varphi_1, \varphi_2 \in \RR_+$.
\ Consequently,
 \[
   \bigl(Z^T_{\Im(\lambda),0}, Z^T_{\Im(\lambda),m_\vartheta} 
         : \lambda \in \Lambda_\vartheta\cap(\ii\RR)\bigr)
   \distr 
   \bigl(\cZ_{\Im(\lambda),0}, \cZ_{\Im(\lambda),m_\vartheta} 
         : \lambda \in \Lambda_\vartheta\cap(\ii\RR)\bigr) , \qquad
   \text{as \ $T \to \infty$,}
 \]
 and hence, by It\^o's formula and the continuous mapping theorem,
 \[
   Y^T \distr \cY \qquad \text{as \ $T \to \infty$}
 \]
 with
 \[
   \cY(t)
   = \sum_{\underset{\SC\tm_\vartheta(\lambda)=m_\vartheta^*}
                    {\lambda \in \Lambda_\vartheta\cap(\ii\RR)}}
      c_{\vartheta,\lambda,m_\vartheta^*}
      \int_0^t
       \cZ_{\Im(\lambda),m_\vartheta^*}(s)
       \, \dd \overline{\cZ_{\Im(\lambda),0}(s)} .
 \]
Further, by Corollary 4.12 in Gushchin and K\"uchler \cite{GusKuc1999},
 \[
   (Y^T(1), [Y^T, Y^T](1)) \distr (\cY(1), [\cY, \cY](1)) \qquad 
   \text{as \ $T \to \infty$.}
 \]
We also have
 \begin{align*}
  \cY(t)
  &= \bbone_{\{0 \in \Lambda_\vartheta,\;\tm_\vartheta(0)=m_\vartheta^*\}}
     c_{\vartheta,0,m_\vartheta^*}
     \int_0^t \cZ_{0,m_\vartheta^*}(s) \, \dd \cW(s) \\
  &\quad
     + \frac{1}{\sqrt{2}}
       \sum_{\underset{\SC\tm_\vartheta(\lambda)=m_\vartheta^*}
                      {\lambda \in \Lambda_\vartheta\cap(\ii\RR_{++})}}
        \biggl[\int_0^t
                \Bigl(c_{\vartheta,\lambda,m_\vartheta^*}
                      \cZ_{\Im(\lambda),m_\vartheta^*}(s)
                      + \overline{c_{\vartheta,\lambda,m_\vartheta^*}
                                  \cZ_{\Im(\lambda),m_\vartheta^*}(s)}\Bigr)
                \, \dd \cW_{\Im(\lambda),\Re}(s) \\
  &\phantom{\quad
     + \frac{1}{\sqrt{2}}
       \sum_{\underset{\SC\tm_\vartheta(\lambda)=m_\vartheta^*}
                      {\lambda \in \Lambda_\vartheta\cap(\ii\RR_{++})}}
        \biggl[}      
               + \int_0^t
                  \Bigl(-c_{\vartheta,\lambda,m_\vartheta^*}
                         \cZ_{\Im(\lambda),m_\vartheta^*}(s)
                        + \overline{c_{\vartheta,\lambda,m_\vartheta^*}
                                    \cZ_{\Im(\lambda),m_\vartheta^*}(s)}\Bigr)
                  \, \dd \cW_{\Im(\lambda),\Im}(s)\biggr] \\
  &= \bbone_{\{0 \in \Lambda_\vartheta,\;\tm_\vartheta(0)=m_\vartheta^*\}}
     c_{\vartheta,0,m_\vartheta^*}
     \int_0^t \cZ_{0,m_\vartheta^*}(s) \, \dd \cW(s) \\
  &\quad
     + \sqrt{2}
       \sum_{\underset{\SC\tm_\vartheta(\lambda)=m_\vartheta^*}
                      {\lambda \in \Lambda_\vartheta\cap(\ii\RR_{++})}}
        \biggl[\int_0^t
                \Re\bigl(c_{\vartheta,\lambda,m_\vartheta^*}
                         \cZ_{\Im(\lambda),m_\vartheta^*}(s)\bigr)
                \, \dd \cW_{\Im(\lambda),\Re}(s) \\
  &\phantom{\quad
     + \sqrt{2}
       \sum_{\underset{\SC\tm_\vartheta(\lambda)=m_\vartheta^*}
                      {\lambda \in \Lambda_\vartheta\cap(\ii\RR_{++})}}
        \biggl[}      
               - \int_0^t
                  \Im\bigl(c_{\vartheta,\lambda,m_\vartheta^*}
                           \cZ_{\Im(\lambda),m_\vartheta^*}(s)\bigr)
                  \, \dd \cW_{\Im(\lambda),\Im}(s)\biggr] ,
 \end{align*}
 hence
 \begin{align*}
  [\cY, \cY](1)
  &= \bbone_{\{0 \in \Lambda_\vartheta,\;\tm_\vartheta(0)=m_\vartheta^*\}}
     c_{\vartheta,0,m_\vartheta^*}^2
     \int_0^1 \cZ_{0,m_\vartheta^*}(s)^2 \, \dd s \\
  &\quad
     + 2 \sum_{\underset{\SC\tm_\vartheta(\lambda)=m_\vartheta^*}
                        {\lambda \in \Lambda_\vartheta\cap(\ii\RR_{++})}}
          \biggl[\int_0^1
                  \Re\bigl(c_{\vartheta,\lambda,m_\vartheta^*}
                           \cZ_{\Im(\lambda),m_\vartheta^*}(s)\bigr)^2
                  \, \dd s
                 + \int_0^1
                    \Im\bigl(c_{\vartheta,\lambda,m_\vartheta^*}
                             \cZ_{\Im(\lambda),m_\vartheta^*}(s)\bigr)^2
                    \, \dd s\biggr] \\
  &= \sum_{\underset{\SC\tm_\vartheta(\lambda)=m_\vartheta^*}
                    {\lambda \in \Lambda_\vartheta\cap(\ii\RR)}}
      |c_{\vartheta,\lambda,m_\vartheta^*}|^2
      \int_0^1 |\cZ_{\Im(\lambda),m_\vartheta^*}(s)|^2 \, \dd s .               
 \end{align*}
Recall that \ $I_3(T) \stoch 0$ \ as \ $T \to \infty$, \ which also implies
 \ $I_1(T) \stoch 0$ \ as \ $T \to \infty$.
\ Further,
 \[
   |I_2(T)|
   \leq \sqrt{\frac{1}{T^{2\ell+3}}
              \int_0^T \Bigl|X^T\Bigl(\frac{t}{T}\Bigr)\Bigr|^2 \dd t
              \int_0^T |\tY(t)|^2 \, \dd t}
   = \sqrt{\int_0^1 |X^T(s)|^2 \, \dd s \;
           \frac{1}{T^2}
           \int_0^T |\tY(t)|^2 \, \dd t}
   \stoch 0
 \]
 as \ $T \to \infty$, \ hence we obtain
 \[
   (\Delta_{\vartheta,T}, J_{\vartheta,T}) \distr (\Delta_\vartheta, J_\vartheta)
   \qquad \text{as \ $T \to \infty$.}
 \]
Moreover, we have \ $J_\vartheta > 0$ \ almost surely.
Indeed, \ $J_\vartheta = 0$ \ would imply
 \ $\int_0^1 |\cZ_{\Im(\lambda),m_\vartheta^*}(s)|^2 \, \dd s = 0$
 \ for all
 \ $\lambda \in \{\lambda \in \Lambda_\vartheta\cap(\ii\RR) : \tm_\vartheta(\lambda)=m_\vartheta^*\}$, 
 \ which, in turn, would imply
 \ $\cZ_{\Im(\lambda),m_\vartheta^*}(s) = 0$ \ for all \ $s \in [0,1]$. 
\ But this is in a contradiction with the fact that
 \ $\cZ_{\Im(\lambda),m_\vartheta^*}$ \ is a non-degenerate Gaussian process.
\proofend

\noindent{\bf Proof of Theorem \ref{PLAMN}.}
We have
 \[
   J_{\vartheta,T}
   = T^{-2m_\vartheta^*} \ee^{-2v_\vartheta^* T}
     \int_0^T Y^{(\vartheta)}(t)^2 \, \dd t ,
   \qquad T \in \RR_+ .
 \]
The process \ $(Y^{(\vartheta)}(t))_{t\in[r,\infty)}$ \ admits the
 representation \eqref{Yvartheta}.
We choose \ $c < v_\vartheta^*$ \ with
 \ $c > \sup\{\Re(\lambda)
              : \lambda \in \Lambda_\vartheta , \;
                \Re(\lambda) < v_\vartheta^*\}$,
\ and apply the representation \eqref{intres}.
By the definition of \ $v_\vartheta^*$, \ we obtain \ $P_{\vartheta,\lambda} = 0$
 \ for each \ $\lambda \in \Lambda_\vartheta$ \ with
 \ $\Re(\lambda) > v_\vartheta^*$, \ hence we obtain
 \begin{equation}\label{yres_PLAMN}
  y_\vartheta(t)
  = \sum_{\lambda\in\Lambda_\vartheta\cap(v_\vartheta^*+\ii\RR)}
     \sum_{\ell=0}^{\tm_\vartheta(v_\vartheta^*)}
      c_{\vartheta,\lambda,\ell} \, t^\ell \ee^{\lambda t} 
     + \Psi_{\vartheta,c}(t) , \qquad t \in \RR_+ .
 \end{equation}
For each \ $\lambda \in \Lambda_\vartheta$ \ and
 \ $\ell \in \{0, \ldots, \tm_\vartheta(\lambda)\}$, \ we have
 \begin{align*}
  c_{\vartheta,\lambda,\ell} \, t^\ell \ee^{\lambda t}
  + c_{\vartheta,\overline{\lambda},\ell} \, t^\ell
    \ee^{\overline{\lambda}t}
  &= t^\ell \bigl(c_{\vartheta,\lambda,\ell} \, \ee^{\lambda t}
                  + \overline{c_{\vartheta,\lambda,\ell}
                    \, \ee^{\lambda t}}\bigr) \\
  &= 2 t^\ell \Re\bigl(c_{\vartheta,\lambda,\ell} \, \ee^{\lambda t}\bigr)
   = t^\ell
     \bigl[\Re\bigl(c_{\vartheta,\lambda,\ell} \, \ee^{\lambda t}\bigr)
           + \Re\bigl(c_{\vartheta,\overline{\lambda},\ell}
                      \, \ee^{\overline{\lambda}t}\bigr)\bigr] ,
 \end{align*}
 hence \eqref{yres_PLAMN} can also be written in the form
 \[
   y_\vartheta(t)
   = \sum_{\lambda\in\Lambda_\vartheta\cap(v_\vartheta^*+\ii\RR)}
      \sum_{\ell=0}^{\tm_\vartheta(\lambda)}
       t^\ell \Re\bigl(c_{\vartheta,\lambda,\ell} \, \ee^{\lambda t}\bigr)
     + \Psi_{\vartheta,c}(t) , \qquad t \in \RR_+ .
 \]
Consequently, by the representation \eqref{Yvartheta}, we have
 \begin{equation}\label{tYres_PLAMN}
  Y^{(\vartheta)}(t) = Y(t) + \tY(t) , \qquad t \in \RR_+ ,
 \end{equation}
 with
 \begin{align*}
  Y(t)
  := \sum_{\lambda\in\Lambda_\vartheta\cap(v_\vartheta^*+\ii\RR)}
      \sum_{\ell=0}^{\tm_\vartheta(\lambda)}
       Y_{\vartheta,\lambda,\ell}(t) , \qquad t \in \RR_+ ,
 \end{align*}
 where the continuous processes
 \ $(Y_{\vartheta,\lambda,\ell}(t))_{t\in\RR_+}$, \ $\lambda \in \Lambda$,
 \ $\ell \in \{0, \ldots, \tm_\vartheta(v_\vartheta^*)\}$, \ and
 \ $(\tY(t))_{t\in\RR_+}$ \ admit representation \eqref{Y} on \ $[r, \infty)$
 \ with \ $y(t) = t^\ell \Re(c_{\vartheta,\lambda,\ell} \ee^{\lambda t})$,
 \ $\lambda \in \Lambda$,
 \ $\ell \in \{0, \ldots, \tm_\vartheta(v_\vartheta^*)\}$, \ and with
 \ $y(t) = \Psi_{\vartheta,c}(t)$, \ $t \in \RR_+$, \ respectively.
The aim of the following discussion is to show that 
 \ $T^{-2m_\vartheta^*} \ee^{-2v_\vartheta^* T} \int_0^T \tY(t)^2 \, \dd t
    \stoch 0$
 \ as \ $T \to \infty$.  
\ We have 
 \[
   \tY(t)
   := \Psi_{\vartheta,c}(t) X_0(0) + \vartheta I_\vartheta(t) + I_\Psi(t) ,
   \qquad t \in [r, \infty) ,
 \]
 with
 \begin{gather*}
  I_\vartheta(t)
  := \int_{[-r,0]} \int_v^0 
      \Psi_{\vartheta,c}(t + v - s) X_0(s) \, \dd s \, a(\dd v) , \qquad
  I_\Psi(t) := \int_0^t \Psi_{\vartheta,c}(t - s) \, \dd W(s) .
 \end{gather*}
The function \ $(\ee^{-ct} \Psi_{\vartheta,c}(t))_{t\in\RR_+}$ \ is bounded, 
 hence \ $|\Psi_{\vartheta,c}(t)| \leq C \ee^{ct}$ \ for all \ $t \in \RR_+$
 \ with \ $C := \sup_{t\in\RR_+} \ee^{-ct} |\Psi_{\vartheta,c}(t)| < \infty$,
 \ hence
 \[
   T^{-2m_\vartheta^*} \ee^{-2v_\vartheta^* T} 
   \int_r^T \Psi_{\vartheta,c}(t)^2 \, \dd t
   \leq \ee^{-2v_\vartheta^* T} \int_0^T C^2 \ee^{2ct} \, \dd t
   = \frac{C^2}{2c} \ee^{-2(v_\vartheta^*-c)T} 
   \to 0 \qquad \text{as \ $T \to \infty$.}
 \]
Next, by Lemma \ref{CS},
 \[
   \int_r^T I_\vartheta(t)^2 \, \dd t
   \leq r \|a\| \int_{-r}^0 X_0(s)^2 \, \dd s
        \int_0^T \psi_\vartheta(v)^2 \, \dd v ,
 \]
 hence we obtain
 \ $T^{-2m_\vartheta^*} \ee^{-2v_\vartheta^* T} \int_0^T I_\vartheta(t)^2 \, \dd t
    \stoch 0$
 \ as \ $T \to \infty$. 
\ Finally, by \eqref{psi},
 \ $T^{-2m_\vartheta^*} \ee^{-2v_\vartheta^* T} \int_0^T I_\Psi(t)^2 \, \dd t
    \stoch 0$
 \ as \ $T \to \infty$, \ and we conclude
 \ $T^{-2m_\vartheta^*} \ee^{-2v_\vartheta^* T} \int_0^T \tY(t)^2 \, \dd t
    \stoch 0$
 \ as \ $T \to \infty$. 
 
Applying Lemma \ref{super}, we obtain
 \begin{align*}
  &T^{-\ell_1-\ell_2} \ee^{-2v_\vartheta^* T}
   \int_0^T
    Y_{\vartheta,\lambda_1,\ell_1}(t)
    Y_{\vartheta,\lambda_2,\ell_2}(t) \, \dd t \\
  &-\int_0^\infty
     \ee^{-2v_\vartheta^* t}
     \Re(c_{\vartheta,\lambda_1,\ell_1} U^{(\vartheta)}_{\lambda_1}
         \ee^{\ii(T-t)\Im(\lambda_1)})
     \Re(c_{\vartheta,\lambda_2,\ell_2} U^{(\vartheta)}_{\lambda_2}
         \ee^{\ii(T-t)\Im(\lambda_2)})
     \dd t
   \as 0
 \end{align*}
 as \ $T \to \infty$ \ for each \ $\lambda_1, \lambda_2 \in \Lambda$ \ with
 \ $\Re(\lambda_1) = \Re(\lambda_2) = v_\vartheta^*$ \ and
 \ $\ell_1 \in \{0, \ldots, \tm_\vartheta(\lambda_1)\}$,
 \ $\ell_2 \in \{0, \ldots, \tm_\vartheta(\lambda_2)\}$.
\ Consequently,
 \[
   T^{-2m_\vartheta^*} \ee^{-2v_\vartheta^* T} \int_0^T Y(t)^2 \, \dd t
   - \int_0^\infty
      \ee^{-2v_\vartheta^* t}
      \Biggl(\sum_{\lambda\in\Lambda_\vartheta\cap(v_\vartheta^*+\ii\RR)}
              \Re\bigl(c_{\vartheta,\lambda,m_\vartheta^*}
                       U^{(\vartheta)}_\lambda
                       \ee^{\ii(T-t)\Im(\lambda)}\bigr)\Biggr)^2
      \dd t
   \as 0
 \]
 as \ $T \to \infty$, \ hence we obtain
 \ $J_{\vartheta,T}- J_\vartheta(T) \as 0$ \ as \ $T \to \infty$.
\ Since \ $H_\vartheta \ne \emptyset$ \ and the numbers in \ $H_\vartheta$ 
 \ have a common divisor \ $D_\vartheta$, \ the process
 \ $(J_\vartheta(t))_{t\in\RR_+}$ \ is periodic with period \ $\frac{2 \pi}{D_\vartheta}$,
 \ and, by Theorem VIII.5.42 of Jacod and Shiryaev \cite{JSh}, we conclude
 \[
   (\Delta_{\vartheta,kD+d}, J_{\vartheta,kD+d})
   \distr (\Delta_\vartheta(d), J_\vartheta(d)) \qquad
   \text{as \ $k \to \infty$}
 \]
 for all \ $d \in \big[0, \frac{2 \pi}{D_\vartheta}\big)$.
\ Moreover, we have \ $J_\vartheta(d) > 0$ \ almost surely for all
 \ $d \in \big[0, \frac{2 \pi}{D_{\vartheta}}\big)$.
\ Indeed, if \ $J_\vartheta(d) = 0$ \ almost surely for all
 \ $d \in \big[0, \frac{2 \pi}{D_{\vartheta}}\big)$, \ then
 \[
   \Re\Biggl(\sum_{\underset{\SC\tm_\vartheta(\lambda)=m_\vartheta^*}
                     {\lambda\in\Lambda_\vartheta\cap(v_\vartheta^*+\ii\RR)}}
                c_{\vartheta,\lambda,m_\vartheta^*}
                U^{(\vartheta)}_\lambda
                \ee^{\ii(d-t)\Im(\lambda)}\Biggr)
   = 0
 \]
  for all \ $d \in \big[0, \frac{2 \pi}{D_{\vartheta}}\big)$ \ and 
  \ $t \in \RR_{++}$.
\ But this is in a contradiction with the fact that the left-hand side is a Gaussian
 random variable with variance  
 \[
   \int_0^\infty
    \Re\left(\sum_{\lambda\in\Lambda_\vartheta\cap(v_\vartheta^*+\ii\RR)}
              c_{\vartheta,\lambda,m_\vartheta^*} \ee^{\ii (d-t) \Im(\lambda)}
              \ee^{-\lambda s} \right)^2 \dd s
   \ne 0 .
 \]

 Consequently we obtain that the family \ $(\cE_T)_{T\in\RR_{++}}$ \ is PLAMN at
 \ $\vartheta$.

If \ $H_\vartheta = \emptyset$, \ then
 \ $\Lambda_\vartheta\cap(v_\vartheta^*+\ii\RR) = \{v_\vartheta^*\}$ \ and
 \ $\tm_\vartheta(v_\vartheta^*) = m_\vartheta^*$, \ thus
 \[
   T^{-2m_\vartheta^*} \ee^{-2v_\vartheta^* T} \int_0^T Y(t)^2 \, \dd t
   \as \int_0^\infty
        \ee^{-2v_\vartheta^* t}
        \bigl(c_{\vartheta,v_\vartheta^*,m_\vartheta^*}
              U^{(\vartheta)}_\lambda\bigr)^2
        \dd t
   = J_\vartheta , \qquad \text{as \ $T \to \infty$.}
 \]
By Theorem VIII.5.42 of Jacod and Shiryaev \cite{JSh}, we conclude
 \[
   (\Delta_{\vartheta,T}, J_{\vartheta,T})
   \distr (\Delta_\vartheta, J_\vartheta) \qquad
   \text{as \ $T \to \infty$,}
 \]
 \ and we obtain that the family \ $(\cE_T)_{T\in\RR_{++}}$ \ is LAMN at
 \ $\vartheta$.
\proofend

\end{document}